\documentclass[french]{amsart}
\usepackage{babel, varioref}
\usepackage{amssymb}
\usepackage[all]{xy}
\newcommand\datver[1]{\def\datverp%
 {\par\boxed{\boxed{\text{#1; Compil\'{e}: \today}}}}}
\numberwithin{equation}{section}
\newtheorem{lemme}{Lemme}[section]
\newtheorem{proposition}[lemme]{Proposition}
\newtheorem{corollaire}[lemme]{Corollaire}
\newtheorem{theorem}{Th\'{e}or\`{e}me}
\newtheorem{non-theorem}[lemme]{Non-Theorem}

\newtheorem{question}[lemme]{Question}

\newtheorem{definition}[lemme]{D\'{e}finition}
\newtheorem{remarque}[lemme]{Remarque}
\newtheorem*{remerciements}{Remerciements}

\newcommand\cf{cf\@. }

\newcommand\eg{e\@.g\@. }
\newcommand\al{al\@.}
\renewcommand\det{\operatorname{det}}
\newcommand\pa{ \partial}

\newcommand\bbB{\mathbb B}
\newcommand\bbC{\mathbb C}

\newcommand\bbN{\mathbb N}
\newcommand\bbP{\mathbb P}
\newcommand\bbQ{\mathbb Q}
\newcommand\bbR{\mathbb R}
\newcommand\bbS{\mathbb S}

\newcommand\bbZ{\mathbb Z}

\newcommand\cU{\mathcal{U}}
\newcommand\cV{\mathcal{V}}

\newcommand\End{\operatorname{End}}

\newcommand\Go{G^{0}}
\newcommand\Goi{G^{0}_{\Id}}
\newcommand\Gos{G^{0}_{\sus}}
\newcommand\Gosu{G^{0}_{\susu}}
\newcommand\Gcs{G^{-1}_{\sus}}
\newcommand\Gc{G^{-1}}
\newcommand\Id{\operatorname{Id}}
\newcommand\CI{\mathcal{C}^{\infty}}
\newcommand\Ld{\operatorname{L}^{2}}
\newcommand\Si{ S_{0}(X)}
\newcommand\Ss{ S_{\sus}(X)}
\newcommand\susn[1]{s(#1)} 
\newcommand\ind{\operatorname{ind}}
\newcommand\tf{\tilde{f}}
\newcommand\hf{\hat{f}}
\newcommand\cF{\mathcal{F}}
\newcommand\cH{\mathcal{H}}
\newcommand\Bo{ B_{0}^{k}}
\newcommand\Bu{ \overline{B}_{1}^{k}}
\newcommand\pt{\operatorname{pt}}
\newcommand\tg{\tilde{g}}
\newcommand\sus{\operatorname{s}(l)}
\newcommand\susu{\operatorname{sus}}
\newcommand\susf{\operatorname{s}}
\newcommand\Hom{\operatorname{Hom}}
\newcommand\detr{\det_{R}}
\newcommand\Trr{\Tr_{R}}
\newcommand\Tr{\operatorname{Tr}}
\newcommand\tr{\operatorname{tr}}
\newcommand\Ar{\operatorname{A_{R}}}
\newcommand\cp[1]{\bbC\bbP_{#1}}  
\newcommand\spec{\operatorname{spec}}
\newcommand\Ch{\operatorname{Ch}}

\renewcommand\Re{\operatorname{Re}}

\datver{0.1D; R\'{e}vis\'{e}: Le 10 mai  2007}
\begin{document}
\title[Op\'{e}rateurs inversibles d'ordre 0]
{Sur la topologie de l'espace des op\'{e}rateurs pseudodiff\'{e}rentiels inversibles d'ordre 0}

\author{Fr\'ed\'eric Rochon}
\address{Department of Mathematics, State University of New York, Stony Brook}
\email{rochon@math.sunysb.edu}
\thanks{L'auteur tient \`a remercier le Fonds qu\'eb\'ecois de la recherche sur la nature et les
technologies pour son soutien financier}
\dedicatory{\datverp}
\begin{abstract}
Les groupes d'homotopie du groupe (stabilis\'e) $\Go(X)$
des op\'erateurs pseudodiff\'erentiels inversibles d'ordre z\'ero
agissant sur une vari\'et\'e compacte sans bord $X$ sont calcul\'es en termes de la $K$-th\'eorie
du fibr\'e cosph\'erique $S^{*}X$.  Du m\^eme coup, on montre que le sous-groupe des perturbations
compactes inversibles de l'identit\'e est faiblement r\'etractile dans $\Go(X)$.  Les r\'esultats
sont aussi adapt\'es au cas des op\'erateurs suspendus.  Des applications \`a la th\'eorie de l'indice et pour le d\'eterminant r\'esiduel de Simon Scott sont aussi donn\'ees.   
\end{abstract}

\maketitle

\tableofcontents

\section*{Introduction}

Soit $X$ une vari\'et\'e compacte et sans bord de classe $\CI$ de dimension sup\'erieure \`a z\'ero.  Le groupe $\Go(X)$ des op\'erateurs
pseudodiff\'erentiels inversibles d'ordre z\'ero agissant sur $\CI(X)$ est un objet important
en g\'eom\'etrie et en analyse.  En th\'eorie de l'indice, la version suspendue
de ce groupe appara\^it lorsqu'on veut d\'ecrire l'op\'erateur normal d'un op\'erateur \`a cusp 
fibr\'e totalement elliptique.  C'est aussi sur le groupe $\Go(X)$ (ou un espace reli\'e)
que plusieurs fonctionnelles jouant le r\^ole de d\'eterminant ont \'et\'e introduites et \'etudi\'ees,
voir par exemple les travaux de Kontsevich et Vishik \cite{Kontsevich-Vishik},
de Scott \cite{Scott}, de Paycha et Scott \cite{Paycha-Scott} et de Friedlander et Guillemin \cite{Friedlander-Guillemin}.  Dans le cas
des op\'erateurs suspendus, Melrose dans \cite{mr96h:58169} a d\'efini
sur le groupe $\Gosu(X)$ des op\'erateurs suspendus inversibles d'ordre
z\'ero une fonctionnelle jouant le r\^ole de l'invariant eta introduit par
Atiyah, Patodi et Singer \cite{APS1}.  Cette fonctionnelle a \'et\'e \'etudi\'ee  entre autres dans les travaux de Moroianu \cite{Moroianu3}, de Lesch, Moscovici et Pflaum \cite{Lesch-Moscovici-Pflaum}
et de Melrose et \al \,  \cite{mr96h:58169}, \cite{Melrose-Nistor},\cite{perdet},\cite{beidb}.

Il appara\^it donc souhaitable d'avoir une bonne compr\'ehension topologique du groupe 
$\Go(X)$.  Dans cet article, on se propose 
d'utiliser les m\'ethodes d\'evelopp\'ees dans \cite{bpfco} pour \'etudier la topologie du
groupe $\Go(X)$ des op\'erateurs pseudodiff\'erentiels inversibles d'ordre z\'ero.  Pour ce faire,
on doit dans un premier temps stabiliser la situation, c'est-\`a-dire permettre \`a ces op\'erateurs
d'agir sur un fibr\'e vectoriel complexe de rang arbitrairement grand.  Dans ce cas, on peut 
calculer les groupes d'homotopie de ce groupe en termes de la $K$-th\'eorie du fibr\'e 
cosph\'erique $S^{*}X$ de $X$ (Th\'eor\`eme~\ref{gh.1}).  Tout comme dans \cite{bpfco}, on
remarque qu'il y a une p\'eriodicit\'e, \`a savoir que les groupes d'homotopie pairs et impairs
sont isomorphes entre eux.  Ce r\'esultat est obtenu en consid\'erant le sous-groupe 
$\Gc(X)\subset \Go(X)$ des perturbations compactes inversibles de l'identit\'e.  Celui-ci,
avec le symbole principal, d\'etermine une fibration de Serre \`a laquelle est associ\'ee une
longue suite exacte de groupes d'homotopie.  En montrant que l'homomorphisme de bord est 
surjectif, on peut alors obtenir le r\'esultat sans trop de difficult\'e.  Notons que dans le cas
o\`u $X=\bbS^{1}$ est un cercle, Melrose dans (\cite{MelroseLMA}, \S~12) a montr\'e qu'il est
possible d'obtenir un sous-groupe faiblement contractile de $\Go(\bbS^{1})$ en imposant certaines
contraintes suppl\'ementaires sur le symbole principal (voir plus bas la discussion dans le \S~\ref{gh.0}).   

Une autre cons\'equence int\'eressante de la surjectivit\'e de l'homomorphisme de bord est 
que le groupe $\Gc(X)$ est faiblement r\'etractile dans $\Go(X)$ (Th\'eor\`eme~\ref{gh.3}), 
c'est-\`a-dire que si $M$ est un $CW$-complexe avec un nombre fini de cellules (\eg une 
vari\'et\'e compacte) et si
$f:M\to \Gc(X)$ est une application continue, alors dans $\Go(X)$, $f$ est homotope
\`a une application constante.   

On montre aussi que la m\'ethode peut \^etre adapt\'ee pour calculer les groupes d'homotopie
de l'espace des op\'erateurs pseudodiff\'erentiels inversibles $l$ fois suspendus.  Dans ce cas,
les groupes d'homotopie sont exprim\'es \`a partir de la $K$-th\'eorie du fibr\'e cosph\'erique
\[
     S^{*}_{X}(X\times \bbR^{l}):= (T^{*}X\times \bbR^{l}\setminus \{0\}) /\bbR^{+}.
\]
Ces r\'esultats donnent lieu \`a une application reli\'ee \`a la th\'eorie de l'indice pour 
les op\'erateurs \`a cusp
fibr\'e.  On montre entre autres qu'un op\'erateur \`a cusp fibr\'e totalement elliptique 
ayant un op\'erateur normal dont le symbole principal est donn\'e par l'identit\'e peut toujours \^etre d\'eform\'e
(apr\`es stabilisation) par une famille d'op\'erateurs totalement elliptiques de sorte que
son op\'erateur normal devienne l'identit\'e, au prix bien entendu de modifier le symbole principal \`a
l'int\'erieur.  Les th\'eor\`emes~\ref{gh.1} et \ref{gh.3} peuvent aussi \^etre utilis\'es pour montrer que le d\'eterminant r\'esiduel
de Simon Scott \cite{Scott} peut \^etre d\'efini globalement sur la composante connexe de l'identit\'e dans $\Go(X)$ (Th\'eor\`eme~\ref{dr.11}).  L'auteur
a r\'ecemment appris que Jean-Marie Lescure et Sylvie Paycha dans
\cite{Lescure-Paycha} ont obtenu un r\'esultat similaire en utilisant des 
m\'ethodes diff\'erentes. 

L'article est organis\'e comme suit.  Dans le \S~\ref{stab.0}, on d\'ecrit comment stabiliser
le groupe $\Go(X)$.  La fibration de Serre qui donne lieu \`a la longue suite exacte de groupes
d'homotopie est ensuite introduite dans le \S~\ref{fdS.0}.  Le \S~\ref{chb.0} donne une description 
de l'homomorphisme de bord en termes d'un indice de famille d'op\'erateurs, ce qui permet de montrer
 sa surjectivit\'e.  Le calcul des groupes d'homotopie de $\Go(X)$ est donn\'e dans le 
 \S~\ref{gh.0}.  Dans le \S~\ref{os.0}, on adapte les r\'esultats au cas des op\'erateurs
 suspendus.  Dans le \S~\ref{ati.0}, on discute d'une application de ces r\'esultats
 en th\'eorie de l'indice pour les op\'erateurs \`a cusp fibr\'e.  Enfin, dans le \S~\ref{dr.0},
 on donne une description topologique du d\'eterminant r\'esiduel de Simon Scott et on montre que ce d\'eterminant admet une d\'efinition globale.

\begin{remerciements}
L'auteur remercie cordialement Richard Melrose pour plusieurs discussions stimulantes sur le sujet.
L'auteur remercie aussi Sergiu Moroianu pour plusieurs remarques et suggestions importantes, notamment concernant le dernier paragraphe de cet article.  Enfin, l'auteur tient \`a remercier le referee pour 
ses commentaires constructifs et d\'etaill\'es.

\end{remerciements}

\section{Stabilisation}\label{stab.0}

Soit $X$ une vari\'{e}t\'{e} compacte sans bord de classe $\CI$ de dimension sup\'erieure \`a z\'ero et soit $E\to X$ un fibr\'{e} vectoriel complexe
sur $X$.  Dans cet article, on se propose d'\'{e}tudier la topologie du groupe
\begin{equation}
  \Go(X;E):= \{ P \in \Psi^{0}(X;E)\quad | \quad P \; \mbox{est inversible} \}
\label{stab.1}\end{equation}
des op\'{e}rateurs pseudodiff\'{e}rentiels (classiques polyhomog\`enes) 
inversibles d'ordre $0$ agissant sur les sections
de $E$.  Plus pr\'{e}cis\'{e}ment, on calculera les groupes d'homotopie de ce groupe muni
de la topologie induite par celle de $\Psi^{0}(X;E)$.  Comme mentionn\'{e} dans l'introduction,
on doit toutefois d'abord stabiliser la situation, c'est-\`{a}-dire permettre au fibr\'{e} vectoriel
$E$ d'avoir un rang arbitrairement large.  De cette fa\c{c}on, la p\'{e}riodicit\'{e} de Bott entre
en jeu et permet d'obtenir un r\'{e}sultat relativement simple.  En fait, en th\'{e}orie de l'indice,
c'est vraiment cette situation qui est d'int\'{e}r\^{e}t.  

Pour r\'{e}aliser une telle stabilisation concr\`{e}tement, soit $\bbS^{1}$ le cercle unit\'{e} dand
$\bbC$. Suivant une id\'{e}e de Richard Melrose (\cf \cite{MelroseLMA}), 
consid\'{e}rons \`{a} la place de $\Go(X;E)$ le groupe d'op\'{e}rateurs
\begin{equation}
\Go(X):= \{ \Id + Q \quad | \quad Q\in \CI(\bbS^{1}\times \bbS^{1}; \Psi^{0}(X)),
\quad \Id +Q \; \mbox{est inversible}\}
\label{stab.2}\end{equation}
agissant sur $\CI(X\times \bbS^{1})$, un op\'{e}rateur $Q\in \CI(\bbS^{1}\times \bbS^{1};
\Psi^{0}(X))$ agissant sur $f\in \CI(X\times \bbS^{1})$ par
\[
   (Qf)(x,\theta)= \int_{\bbS^{1}} (Q(\theta,\theta') f_{\theta'})(x) d\theta'
\]
o\`u $\theta\in [0,2\pi)$ est la coordonn\'ee angulaire usuelle sur $\bbS^{1}$ et la fonction
$f_{\theta'}\in \CI(X)$ est donn\'{e}e par
\[
            f_{\theta'}(x):= f(x, e^{i\theta'}), \quad x\in X, \quad \theta'\in [0,2\pi),
            \; e^{i\theta'}\in \bbS^{1}.
\]
En termes de la base orthonormale de $\Ld(\bbS^{1})$ donn\'{e}e par les fonctions
propres $\frac{e^{i k\theta}}{\sqrt{2\pi}}$, $k\in \bbZ$ de l'op\'{e}rateur
$\frac{d}{d\theta}$, un op\'erateur $Q\in \CI(\bbS^{1}\times \bbS^{1};\Psi^{0}(X))$
peut \^etre d\'ecrit par une matrice $\bbZ\times \bbZ$ avec coefficients
\begin{equation}
    a_{kl}:= \frac{1}{2\pi} \langle e^{-ik\theta}; Qe^{il\theta}\rangle_{\Ld(\bbS^{1})}
              \in \Psi^{0}(X), \quad k,l\in \bbZ
\label{stab.3}\end{equation}
d\'ecroissant rapidement vers z\'ero lorsque $|k|+|l|\to \infty$.

Soit $F\to X$ un autre fibr\'e vectoriel complexe tel que $E\oplus F$ 
s'identifie avec le fibr\'{e} trivial $\underline{\bbC}^{n}$ de rang $n$ (pour $n$ assez grand).
En identifiant $\bbC^{n}$ avec le sous-espace vectoriel de $\Ld(\bbS^{1})$ ayant pour base
$\{1, e^{i\theta},\ldots, e^{i(n-1)\theta}\}$, on peut alors regarder $E$ comme un sous-fibr\'e
vectoriel de $X\times \Ld(\bbS^{1})\to X$, ce qui donne lieu a un plongement\footnote{
Ce n'est pas n'importe quelle identification de $E$ avec un sous-fibr\'e vectoriel
de $X\times \Ld(\bbS^{1})\to X$ qui donne lieu \`a un tel plongement.}
\begin{equation}
    \begin{array}{ccc} \Go(X;E) & \subset & \Go(X) \\
                            P & \mapsto & \Id + (P-\Id_{E}),
    \end{array}                         
\label{stab.4}\end{equation}
l'op\'erateur $\Id_{E}:X\times \Ld(\bbS^{1})\to X\times\Ld(\bbS^{1})$
\'etant donn\'e par la projection orthogonale sur $E$.

C'est un tel plongement qui permet de voir $\Go(X)$ comme la stabilisation de 
$\Go(X;E)$, puisque dans $\Go(X)$, un op\'erateur peut agir sur un sous-fibr\'e vectoriel
de $X\times \Ld(\bbS^{1})\to X$ de rang arbitrairement large.  D'un autre c\^ot\'e, la condition
de d\'ecroissance rapide sur les coefficients \eqref{stab.3} permet toujours d'approximer 
un op\'erateur de $\Go(X)$ par un autre op\'erateur de $\Go(X)$ agissant sur un sous-fibr\'{e}
vectoriel de $\Ld(\bbS^{1})$ de rang fini.  

Remarquons que pour d\'efinir le groupe $\Go(X)$, on aurait pu tout aussi bien prendre une
vari\'et\'e compacte sans bord \`a la place du cercle.  Cela aurait donn\'e lieu \`a la m\^eme
structure de groupe topologique.  \`A la place de $G^{0}(X)$, on aurait pu aussi prendre la limite
t\'elescopique\footnote{suivant la terminologie de Bott et Tu (\cite{Bott-Tu}, p.241)}
\[
                   \lim_{n\to +\infty} \Go(X; \underline{\bbC}^{n})
\]
d\'efinie via les inclusions $\Go(X;\underline{\bbC}^{n})\subset \Go(X;\underline{\bbC}^{n+1})$ pour
$n\in \bbN$.  \`A strictement parler, la topologie de cet espace est l\'eg\`erement diff\'erente de 
celle de $G^{0}(X)$, mais conduit au m\^eme type d'homotopie.  

\section{La fibration de Serre associ\'ee}\label{fdS.0}

Un des sous-groupes importants de $\Go(X)$ est donn\'e par 
\begin{equation}
  \Gc(X):= \{ \Id +Q \quad | \quad 
           Q\in \CI(\bbS^{1}\times \bbS^{1}; \Psi^{-1}(X)), \;
           \Id +Q \, \mbox{est inversible} \},  
\label{fdS.1}\end{equation}
le sous-groupe des perturbations compactes inversibles de l'identit\'{e}.  Il donne 
lieu \`{a} une suite exacte \`a gauche
\begin{equation}
  0\to \Gc(X) \hookrightarrow \Go(X) \overset{\sigma}{\longrightarrow}
                 \CI(S^{*}X; G^{-\infty}(\bbS^{1}))
\label{fdS.2}\end{equation} 
o\`u $S^{*}X:= (T^{*}X\setminus X)/\bbR^{+}$ est le \emph{fibr\'e cosph\'erique},
\begin{equation*}
   G^{-\infty}(\bbS^{1}):= \{\Id + A \quad | \quad 
                 A\in \Psi^{-\infty}(\bbS^{1}), \; \Id+A\, \mbox{est inversible}\}
\end{equation*}
est le groupe des perturbations r\'egularisantes inversibles de l'identit\'e sur 
$\bbS^{1}$ et 
\[
             \sigma: \Go(X)\to \CI(S^{*}X; G^{-\infty}(\bbS^{1}))
\]
est l'application qui, \`a un op\'erateur donn\'e, lui associe son symbole principal.
L'application $\sigma$ n'est toutefois pa surjective du fait que la condition 
d'inversibilit\'e impose certaines restrictions sur les valeurs possibles du symbole
principal.  En effet, comme le groupe $G^{-\infty}(\bbS^{1})$ est un espace classifiant
pour la $K$-th\'eorie impaire, on obtient une application
\begin{equation*}
    h: \CI(S^{*}X; G^{-\infty}(\bbS^{1}))\to K^{1}(S^{*}X)
\label{fdS.3}\end{equation*}
qui \`a un \'el\'ement $s\in \CI(S^{*}X; G^{-\infty}(\bbS^{1}))$ associe sa classe 
d'homotopie.  Si d'autre part
\[
        \delta: K^{1}(S^{*}X)\to K^{0}_{c}(T^{*}X)\cong K^{0}(
        \overline{T^{*}X}, S^{*}X)
\]
d\'enote l'homomorphisme de bord de la suite exacte \`a six termes associ\'ee 
\`a la paire $(\overline{T^{*}X},S^{*}X)$ o\`u 
$S^{*}X \subset \overline{T^{*}X}$ est vu comme le bord de la 
\emph{compactification radiale} $\overline{T^{*}X}$ du fibr\'e cotangent
$T^{*}X$, alors l'indice d'un op\'erateur $\Id+Q\in \Go(X)$, via le th\'eor\`eme 
d'Atiyah-Singer \cite{Atiyah-Singer1} est donn\'e par
\[
    \ind_{a}(\Id +Q) = \ind_{t} \circ \delta \circ h\circ\sigma(\Id+Q)
\]
o\`u $\ind_{t}: K^{0}_{c}(T^{*}X)\to \bbZ$ est l'indice topologique de 
Atiyah-Singer.  Comme $\Id+Q$ est par d\'efinition un op\'{e}rateur 
inversible, son indice est n\'ecessairement z\'ero, ce qui impose une restriction sur la 
classe d'homotopie de son symbole principal.  

D\'esignons par $\Si$ le noyau de l'application 
\[
                     \ind_{t}\circ \delta \circ h: \CI(S^{*}X;G^{-\infty}(\bbS^{1}))\to \bbZ.
\]
\begin{lemme}
L'application $\sigma$ associ\'ee au symbole principal d'un op\'erateur donne lieu \`a la
suite exacte
\[
 0\to \Gc(X) \hookrightarrow \Go(X) \overset{\sigma}{\longrightarrow}
              \Si\to 0.
\]
\label{fdS.4}\end{lemme}
\begin{proof}
Par la discussion pr\'ec\'edente, on sait que 
\[
      \sigma(\Go(X))\subset \Si
\]
et que la suite est exacte \`a gauche.  Il suffit donc de montrer que l'application
$\sigma$ est surjective.  Soit $p\in \Si$ un \'el\'ement donn\'e, alors on sait au 
moins qu'il existe $P\in \Id+ \CI(\bbS^{1}\times \bbS^{1}; \Psi^{0}(X))$ avec symbole
principal donn\'e par p:
\[
                   \sigma(P)=p.
\]
Cet op\'erateur n'a a priori aucune raison d'\^etre inversible.  Toutefois, son indice
est n\'ecessairement nul, puisque d'apr\`es le th\'eor\`eme d'Atiyah-Singer 
\cite{Atiyah-Singer1}, celui-ci est donn\'{e} par
\[
            \ind_{a}= \ind_{t}\circ \delta \circ h(p)=0.
\]
Le noyau et le conoyau de $P$ ont donc la m\^eme dimension.
Si $Q: \ker P\to \ker P^{*}$ repr\'esente un choix d'isomorphisme entre ces derniers,
on peut alors, en posant $Q=0$ sur le compl\'ement orthogonal de $\ker P$ dans
$\Ld(X\times \bbS^{1})$, interpr\'eter $Q$ comme \'etant un \'el\'ement de 
$\CI(\bbS^{1}\times \bbS^{1}; \Psi^{-\infty}(X))$.  De cette fa\c{c}on, on obtient que
$P+Q\in \Go(X)$ est inversible avec symbole principal $\sigma(P+Q)=p$, ce qui donne
le r\'esultat cherch\'e. 
\end{proof}

La suite exacte du lemme pr\'{e}c\'edent est en fait une fibration localement triviale.  
Il suffit de choisir une application de quantification
\[
                    q: \CI(S^{*}X; \Id+ \Psi^{-\infty}(\bbS^{1}))\to \Psi^{0}(X).
\]
Comme $\Go(X)$ est ouvert dans $\Psi^{0}(X)$, on obtient alors facilement des sections locales de
$\sigma: \Go(X)\to \Si$ en utilisant l'application de quantification $q$.  En particulier, cette 
fibration localement triviale est une fibration de Serre, ce qui veut dire qu'on a
une longue suite exacte  pour les groupes d'homotopie
\begin{multline}
\cdots \pi_{k}(\Gc(X))\to \pi_{k}(\Go(X))\overset{\sigma}{\to} \pi_{k}(\Si)
   \overset{\pa}{\to} \pi_{k-1}(\Gc(X)) \to \cdots \\
 \cdots \to \pi_{1}(\Si) \overset{\pa}{\to} \pi_{0}(\Gc(X))\to \pi_{0}(\Go(X))
 \overset{\sigma}{\to} \pi_{0}(\Si).
\label{fdS.5}\end{multline}
C'est par le biais de cette longue suite exacte que nous allons calculer les groupes d'homotopie
de $\Go(X)$.  En effet, on sait que $\Gc(X)$ est un espace classifiant pour la $K$-th\'eorie 
impaire, donc ses groupes d'homotopie sont donn\'{e}s par
\begin{equation}
\pi_{k}(\Gc(X))\cong \left\{ 
              \begin{array}{cl}
                 \bbZ, & k \quad \mbox{impair}, \\
                 \{0\},& k \quad \mbox{pair}.
              \end{array} \right.
\label{fdS.6}\end{equation}
De plus, comme $G^{-\infty}(\bbS^{1})$ est aussi un espace classifiant pour la $K$-th\'eorie
impaire, on peut v\'erifier que pour $k>0$,
\begin{equation}
\begin{aligned}
   \pi_{k}(\Si)&\cong [ S^{k}(S^{*}X); G^{-\infty}(\bbS^{1})] \\
               &\cong K^{-k-1}(S^{*}X)
\end{aligned}
\label{fdS.7}\end{equation} 
o\`u $S^{k}(S^{*}X)$ d\'{e}note la $k$-suspension de $S^{*}X$, alors que pour $k=0$, on a
\begin{equation}
  \pi_{0}(\Si)\cong \ker [\ind_{t}\circ \delta: K^{1}(S^{*}X)\to \bbZ]
\label{fdS.8}\end{equation} 
essentiellement par d\'efinition de $\Si$.  

Une bonne compr\'ehension de l'homomorphisme de bord
\[
          \pa: \pi_{k}(\Si)\to \pi_{k-1}(\Gc(X))
\]
nous permettra donc de calculer les groupes d'homotopie de $\Go(X)$.

\section{Caract\'erisation de l'homomorphisme de bord $\pa$}\label{chb.0}

Dans cette section, nous allons interpr\'eter l'homomorphisme de bord comme \'etant
un indice de famille d'op\'erateurs de Fredholm.  \'Etant donn\'e une application 
$f: \bbS^{k}\to \Si$ envoyant le point de base de $\bbS^{k}$ (choisi au pr\'ealable) sur l'application
identit\'e, on peut relever celle-ci dans l'espace $\Id +\CI(\bbS^{1}\times \bbS^{1};\Psi^{0}(X))$, 
c'est-\`a-dire qu'il existe une application 
\[
     \tf: \bbS^{k}\to \Id +\CI(\bbS^{1}\times \bbS^{1}; \Psi^{0}(X)) 
\]  
telle que $\sigma\circ \tf= f$.  Comme le symbole principal de la famille d'op\'erateurs 
d\'efinie par $\tf$ est inversible, l'application $\tf$ d\'efinit une famille d'op\'erateurs
de Fredholm agissant sur l'espace de Hilbert $\cH:= \Ld(X\times \bbS^{1})$ que l'on d\'enotera
\[
    \hf: \bbS^{k}\to \cF(\cH)
\]
o\`u $\cF(\cH)$ repr\'esente l'espace des op\'erateurs (born\'es) de Fredholm agissant sur
$\Ld(X\times \bbS^{1})$.  En choisissant $\tf$ de fa\c{c}on appropri\'ee, on peut toujours supposer
que $\tf$ envoie le point de base de $\bbS^{k}$ sur l'identit\'e dans $\cF(\cH)$.  Comme on peut
le v\'erifier ais\'ement, la classe d'homotopie $[\hf]\in \pi_{k}(\cF(\cH))$ d\'efinie par
$\hf$ ne d\'epend que de la classe d'homotopie  $[f]\in \pi_{k}(\Si)$ associ\'ee \`a $f$.  
De cette fa\c{c}on, on d\'efinit donc une application
\begin{equation}
   \begin{array}{lccc}
          q: & \pi_{k}(\Si) & \to & \pi_{k}(\cF(\cH)) \\
             &       [f] & \mapsto  & [\hf].
   
   \end{array}
\label{chb.1}\end{equation}
D'un autre c\^ot\'e, comme $\Gc(X)$ est un espace classifiant pour la $K$-th\'eorie 
impaire, on a une identification
\begin{equation}
\pi_{k-1}(\Gc(X))\cong \widetilde{K}^{-1}(\bbS^{k-1})\cong \widetilde{K}^{0}(\bbS^{k})
\label{chb.2}\end{equation}
o\`u $\widetilde{K}$ d\'enote la $K$-th\'eorie r\'eduite pour un espace muni d'un point
de base.  
\begin{proposition}
Vu comme une application
\[
     \pa: \pi_{k}(\Si) \to \widetilde{K}^{0}(\bbS^{k})
\]
en utilisant l'identification~\eqref{chb.2}, l'homomorphisme de bord $\pa$ est donn\'e
par $\pa= \ind\circ q$ o\`u 
\[
     \ind: \pi_{k}(\cF(\cH))\to \widetilde{K}^{0}(\bbS^{k})
\]
est l'indice de famille d'op\'erateurs de Fredholm tel que d\'efini dans l'appendice de 
\cite{Atiyah}.
\label{chb.3}\end{proposition}
\begin{proof}
La preuve est essentiellement la m\^eme que dans (\cite{bpfco}, proposition 8.15).  Nous donnerons tout de m\^eme une preuve en r\'ef\'erant \`a \cite{bpfco} pour plus de d\'etails.  

Soit $f:\bbS^{k}\to \Si$ une application repr\'esentant la classe d'homotopie 
$[f]\in \pi_{k}(\Si)$ et soit $\tf: \Id+ \CI(\bbS^{1}\times \bbS^{1}; \Psi^{0}(X))$ un
choix de rel\`evement, de sorte que vu comme une famille d'op\'erateurs de Fredholm
\[
    \hf: \bbS^{k}\to \cF(\cH),
\]
on ait $q([f])= [\hf]$.  Sans changer la classe d'homotopie de $f$ et $\tf$, on peut supposer que $f\equiv \Id$, $\tf\equiv \Id$ dans une boule ouverte $\Bo\subset \bbS^{k}$ contenant le point de base de $\bbS^{k}$.  Soit $\Bu\subset \bbS^{k}$ le compl\'ement de $\Bo$ dans $\bbS^{k}$. 
Intuitivement, le r\'esultat n'est pas surprenant puisqu'\`a la fois l'homomorphisme de bord $\pa$ et
l'indice de famille $\ind\circ q$ mesurent l'obstruction \`a relever une application $f$
dans $\Go(X)$.  

En prenant un sous-espace $V\subset \mathcal{H}$ de codimension finie tel que 
$\ker \hat{f}(s) \cap V= \{0\}$ pour tout $s\in \bbS^{k}$, on peut d\'ecrire l'indice de 
famille de $\hat{f}$ par 
\[
        \ind(\hat{f})= [V^{\perp}] - [ \hat{f}(V)^{\perp}]
\]                    
o\`u $V^{\perp}$ et $\hat{f}(V)^{\perp}$ sont des fibr\'es vectoriels sur $\bbS^{k}$.
Comme l'espace $\Bu$ est contractile et que $\left. \tf\right|_{\pa \Bu}\equiv \Id$, l'indice de la 
famille d'op\'{e}rateurs $\hf$ restreinte \`a $\Bu$ est nul.  En choisissant $V$ minutieusement, on 
peut aussi supposer que $V^{\perp}$ et $\hf(V)^{\perp}$ sont isomorphes en tant que fibr\'es vectoriels
lorsque restreints \`a $\Bu$ (voir \cite{bpfco}, lemme 8.14).  Soit $\varphi: V^{\perp}\to \hf(V)^{\perp}$ un choix
d'isomorphisme explicite sur $\Bu$.  En posant que $\varphi$ agit par z\'ero sur $V$, on obtient une
famille d'op\'{e}rateurs
\[
            \phi: \Bu\to \CI(\bbS^{1}\times \bbS^{1};\Psi^{-\infty}(X))\cong 
            \Psi^{-\infty}(X\times \bbS^{1}).
\]
Puisque $\Bu$ est compact, il existe $\lambda>0$ tel que
\[
     \hf(s)+ \lambda \phi(s)\in \Go(X), \quad \forall s\in \Bu.
\]
En renormalisant $\phi$ si n\'ecessaire, on peut donc supposer que $\tf+ \phi$ est une
application de la forme
\[
                \tf+\phi: \Bu\to \Go(X).
\]
Puisque $\phi\in \CI(\bbS^{1}\times \bbS^{1};\Psi^{-\infty}(X))\subset 
\CI(\bbS^{1}\times \bbS^{1};\Psi^{-1}(X))$, cette application est aussi un rel\`evement de 
$f$
\[
         \sigma(\tf+\phi)=f \quad \mbox{sur} \; \Bu.
\]
De plus, \'etant donn\'ee que $\tf\cong \Id$ sur $\pa \Bu$, l'application 
$\tf+\phi$ prend valeur dans $\Gc(X)$ lorsque restreinte \`a $\pa\Bu$.  D'autre part,
par d\'efinition de l'homomorphisme de bord (voir \cite{Steenrod}, \S17.1), on a
\[
   \pa([f])= [\left.(\tf+\phi)\right|_{\pa \Bu}] \in \pi_{k-1}(\Gc(X)).
\]
\'Etant donn\'e que $V^{\perp}= \tf(V)^{\perp}$ canoniquement sur $\pa \Bu$, l'application
$\phi$ prend la forme
\[
           \phi: \pa \Bu\to \End( V^{\perp}, V^{\perp}).
\]
Or, clairement, la construction de recollement\footnote{clutching construction en anglais, voir
\cite{Atiyah} pour une description.} appliqu\'ee \`a $(\Id+\phi)^{-1}$ donne le fibr\'e vectoriel
virtuel $[\hf(V)^{\perp}]- [V^{\perp}]$.   La construction de recollement appliqu\'ee \`a 
$(\Id+\phi)$ donne donc $[V^{\perp}]- [T(V)^{\perp}]$ ce qui montre que 
\[
                  \pa= \ind\circ q
\]
puisque l'identification~\eqref{chb.2} est donn\'ee par la construction de recollement.
\end{proof}
\begin{lemme}
L'homomorphisme de bord $\pa: \pi_{k}(\Si)\to \widetilde{K}^{0}(\bbS^{k})$ est
surjectif pour tout $k\in \bbN$.
\label{chb.4}\end{lemme}
\begin{proof}
Par la proposition pr\'ec\'edente, il suffit de montrer que l'application 
$\ind\circ q$ est surjective.  Montrons d'abord que, dans ce qui 
 correspond en quelque sorte au cas $k=0$,
\begin{equation}
   \ind\circ q: \pi_{0}(S^{*}X; G^{-\infty}(\bbS^{1}))\to \bbZ\cong K^{0}(\pt),
\label{chb.5}\end{equation} 
on a aussi une application surjective.  Supposons que $l\in \bbZ$ est donn\'e.  Soit
alors
\[
         g:\bbS^{2n-1}\to G^{-\infty}(\bbS^{1}), \quad n=\dim X,
\]  
une application qui, via la construction de recollement de la s\'erie d'identifications
\begin{equation}
\begin{aligned}
  \pi_{2n-1}(G^{-\infty}(\bbS^{1})) &\cong \widetilde{K}^{0}(\bbS^{2n}), \quad
    \mbox{(construction de recollement)}  \\
    &\cong K^{0}(\pt), \quad \mbox{(P\'eriodicit\'e de Bott)}  \\
    &\cong \bbZ
\end{aligned}
\label{chb.6}\end{equation}
correspond \`a l'entier $l$.  En regardant la sph\`ere $\bbS^{2n-1}$ comme \'etant donn\'ee
par 
\[
     \bbS^{2n-1}\cong \bbB^{2n-1}/ \pa \bbB^{2n-1}
\]
o\`u $\bbB^{2n-1}\subset \bbR^{2n-1}$ est la boule ferm\'ee de rayon $1$, on obtient une fonction
\[
    \tg: \bbB^{2n-1}\to G^{-\infty}(\bbS^{1})
\]
qui envoie le bord de $\bbB^{2n-1}$ sur l'identit\'e.  D'autre part, soit 
\[
     i: \bbB^{2n-1} \hookrightarrow S^{*}X
\]
un plongement de la boule $\bbB^{2n-1}$ dans $S^{*}X$.  L'application
$\tg$ d\'efinit alors une application
\[
   \tg\circ i^{-1}: i(\bbB^{2n-1})\to G^{-\infty}(\bbS^{1})
\]
qui peut \^etre \'etendue \`a tout $S^{*}X$ par l'identit\'e.  Soit 
\[
        f: S^{*}X\to G^{-\infty}(\bbS^{1})
\]
cette extension de $\tg\circ i^{-1}$ \`a tout $S^{*}X$.  Il est alors ais\'e de montrer que
l'indice topologique 
\[
    \ind_{t}\circ \delta([f])
\]
associ\'e \`a la $K$-classe $[f]\in K^{1}(S^{*}X)$ de $f$ est exactement $l$, ce qui 
d\'emontre la surjectivit\'e de \eqref{chb.5}.  Fort de ce r\'esultat, le lemme se 
d\'emontre essentiellement comme dans (\cite{Atiyah}, proposition A6) en utilisant
un op\'erateur
\[
              T\in \Psi^{0}(X;\bbC^{N}) \quad (N\in \bbN \;\mbox{assez grand})
\]
elliptique (donc de Fredholm) d'indice $-1$.  Un tel op\'erateur existe par la
surjectivit\'e de \eqref{chb.5} que nous venons d'\'etablir.
\end{proof}

\section{Les groupes d'homotopie de $\Go(X)$}\label{gh.0}

La surjectivit\'e de l'homomorphisme de bord $\pa$ nous permet maintenant de calculer les groupes
d'homotopie de $\Go(X)$.

\begin{theorem}
Soit $X$ une vari\'et\'e compacte sans bord de dimension sup\'erieure \`a z\'ero.
Les groupes d'homotopie de $\Go(X)$ sont donn\'es par
\[
    \pi_{k}(\Go(X))\cong \left\{ 
      \begin{array}{l}
      \ker[ \ind_{t}\circ \delta: K^{1}(S^{*}X)\to \bbZ ], \quad k \,\mbox{pair}, \\
      K^{0}(S^{*}X), \quad k \,\mbox{impair},
      \end{array}\right.
\]
o\`u $\delta: K^{1}(S^{*}X)\to K^{0}_{c}(T^{*}X)$ est l'homomorphisme de bord et 
$\ind_{t}: K^{0}_{c}(T^{*}X)\to \bbZ$ est l'indice topologique d'Atiyah-Singer. 
\label{gh.1}\end{theorem}
\begin{proof}
Par le lemme~\ref{chb.4}, l'homomorphisme de bord
\[
   \pa: \pi_{k}(\Si)\to \pi_{k-1}(\Gc(X))
\]
est surjectif pour tout $k\in \bbN$.  La longue suite exacte~\ref{fdS.5} se d\'ecompose
donc en courtes suites exactes
\begin{equation}
\begin{gathered}
0\to \pi_{k}(\Go(X))\to \pi_{k}(\Si)\to 0, \quad k \,\mbox{impair}, \\
0\to \pi_{k}(\Go(X))\to \pi_{k}(\Si)\to \bbZ \to 0, \quad k \, \mbox{pair},
\end{gathered}
\label{gh.2}\end{equation}
o\`u on a utilis\'e le fait que $\Gc(X)$ est un espace classifiant pour la $K$-th\'eorie impaire,
\[
     \pi_{k}(\Gc(X))\cong \left\{ 
       \begin{array}{ll}
       \bbZ, \quad k\,\mbox{impair}, \\
       \{0\}, \quad k\, \mbox{pair}.
       \end{array} \right.
\]
Le r\'esultat suit en utilisant l'identification~\eqref{fdS.8}.
Dans le cas o\`u $k=0$, la surjectivit\'e \`a droite est une cons\'equence du lemme~\ref{fdS.4}
et le r\'esultat est alors obtenu en utilisant l'identification~\eqref{fdS.8}.
\end{proof}
\begin{remarque}
Lorsque la vari\'et\'e $X$ est constitu\'ee d'un nombre fini de points, le th\'eor\`eme ne s'applique pas, mais on peut identifier $\Go(X)$ 
avec $G^{-\infty}(\bbS^{1})$, qui est un espace classifiant pour la $K$-th\'eorie impaire.  Les groupes 
d'homotopie pairs sont donc triviaux alors que les groupes d'homotopie impairs sont isomorphes \`a $\bbZ$
dans ce cas.
\label{rem.1}\end{remarque}

Consid\'erons le cas particulier o\`u $X= \bbS^{1}$ est donn\'e par le cercle.  Alors
le fibr\'e cosph\'erique
\[
       S^{*}\bbS^{1}\cong \bbS^{1}_{+}\sqcup \bbS_{-}^{1}
\]
est l'union disjointe de deux cercles.  Par suite, la $K$-th\'eorie de cet espace est
donn\'ee par
\begin{equation*}
\begin{aligned}
K^{0}(S^{*}\bbS^{1})&\cong \bbZ\oplus \bbZ  \\
K^{1}(S^{*}\bbS^{1})&\cong  \bbZ\oplus \bbZ,
\end{aligned} 
\end{equation*}
ce qui donne pour les groupes d'homotopie de $\Go(\bbS^{1})$
\[
\pi_{k}(\Go(\bbS^{1}))\cong \left\{ 
       \begin{array}{ll}
       \bbZ, \quad k\,\mbox{pair}, \\
       \bbZ \oplus \bbZ, \quad k\, \mbox{impair}.
       \end{array} \right.
\]
Lorsque $X=\bbS^{2}$ est la sph\`ere de dimension deux, on peut aussi calculer
explicitement les groupes d'homotopies de $\Go(\bbS^{2})$.  On consid\`ere
d'abord la suite exacte \`a six termes associ\'ee \`a la paire
$(\overline{T^{*}\bbS^{2}}, S^{*}\bbS^{2})$ o\`u $S^{*}\bbS^{2}$ est vu comme
le bord de la compactification radial $\overline{T^{*}\bbS^{2}}$ du fibr\'e
cotangent $T^{*}\bbS^{2}$.  Apr\`es les identifications \'evidentes
\begin{equation}
\begin{gathered}
   K^{i}(\overline{T^{*}\bbS^{2}}, S^{*}\bbS^{2})\cong K^{i}_{c}(T^{*}\bbS^{2}), \\
   K^{i}(\overline{T^{*}\bbS^{2}})\cong K^{i}(\bbS^{2}),
\end{gathered}
\label{hsd.2}\end{equation} 
cette suite exacte prend la forme
\begin{equation}
\xymatrix{K_{c}^{0}(T^{*}\bbS^{2})\ar[r]&
 K^{0}(\bbS^{2})\ar[r]^{\pi^{*}}&
 K^{0}(S^{*}\bbS^{2})\ar[d]\\
 K^{1}(S^{*}\bbS^{2})\ar[u]&
 K^1(S^{2})\ar[l]&
 K_{c}^{1}(T^{*}\bbS^{2})\ar[l]}
\label{hsd.1}\end{equation}
o\`u $\pi: S^{*}\bbS^{2}\to \bbS^{2}$ est la projection de fibr\'e.  Par 
l'isomorphisme de Thom en $K$-th\'eorie, on a que
\begin{equation}
  K^{j}_{c}(T^{*}\bbS^{2}) \cong K^{j}(\bbS^{2}) \cong 
  \left\{ \begin{array}{ll}
            \bbZ\oplus \bbZ, & j=0,\\
            \{0\},   & j=1.
          \end{array} \right.
\label{hsd.3}\end{equation}
D'autre part, en regardant $\bbS^{2}$ comme \'etant $\cp{1}$, on a que
\[
        K^{0}(\bbS^{2}) \cong \frac{\bbZ[t]}{(t-1)^{2}},
\]
l'isomorphisme \'etant donn\'e par $[\underline{\bbC}]\mapsto 1$ et
$[H]\to t$ o\`u $H\to \cp{1}$ est le fibr\'e en droite d'hyperplans.  Clairement,
pour tout $n\in \bbN$, $\pi^{*}[\underline{\bbC^{n}}]$ n'est pas nul dans
$K^{0}(S^{*}\bbS^{2})$.  En identifiant $T^{*}\bbS^{2}$ avec
$H\otimes_{\bbC} H$ et en choisissant une m\'etrique hermitienne sur $H$ (et donc
sur $H\otimes_{\bbC} H$), on peut alors identifier $S^{*}\bbS^{2}$ avec 
$S(H\otimes_{\bbC} H)$, le fibr\'e en cerle unitaire de $H\otimes_{\bbC} H$
\[
      S(H\otimes_{\bbC} H):= \{ v\in H\otimes_{\bbC} H \quad | \quad |v|=1\}.
\]
En ce cas, il devient \'evident que le fibr\'e en droite 
$\pi^{*}(H\otimes_{\bbC} H)\to S^{*}\bbS^{2}$ est trivial sur $S^{*}\bbS^{2}$.  Ainsi,
comme $([H]-1)^{2}=0$, on a que
\[
        \pi^{*}(2[H])= \pi^{*}([H]^{2}-1)=0.
\]
Cependant, $\pi^{*}[H]$ n'est pas trivial puisque d'apr\`es la 
suite de Gysin associ\'ee au fibr\'e $T^{*}\bbS^{2}\to \bbS^{2}$,
sa classe de Chern est le g\'en\'erateur de $H^{2}(S^{*}S^{2})\cong \bbZ_{2}$. 
La suite exacte \eqref{hsd.1}
et l'isomorphisme de Thom \eqref{hsd.3} montrent donc que
\[
  K^{0}(S^{*}\bbS^{2}) \cong \bbZ\oplus\bbZ^{2}, \quad K^{1}(S^{*}(\bbS^{2}))\cong \bbZ,
\]
et donc que les groupes d'homotopie de $\Go(\bbS^{2})$ sont donn\'es
par
\[
     \pi_{k}(\Go(\bbS^{2}))\cong \left\{
          \begin{array}{ll}
               \{0\}, & k \; \mbox{pair},  \\
               \bbZ\oplus\bbZ_{2} & k\; \mbox{impair}.
          \end{array}   \right.
\] 

Plus g\'en\'eralement, en utilisant l'isomorphisme donn\'e par le caract\`ere de Chern
\[
    \Ch: K^{*}(S^{*}X)\otimes_{\bbZ}\bbQ \to H^{*}(S^{*}X,\bbQ),
\]
on peut exprimer les groupes d'homotopie rationnels $\pi_{i}(\Go(X))\otimes_{\bbZ}\bbQ$ de
$\Go(X)$ en termes de la cohomologie paire et impaire de $S^{*}X$.
Notons aussi que comme la $K$-th\'eorie paire (non-r\'eduite) 
n'est jamais triviale, on a que pour toute vari\'et\'e compacte $X$ sans bord, 
le groupe fondamental de 
$\Go(X)$ n'est jamais trivial.  En particulier, le groupe $G^{0}(X)$ n'est 
jamais contractile. Toutefois, comme il est montr\'e par
Melrose (\cite{MelroseLMA}, \S 12), lorsque $X=\bbS^{1}$ est le cercle,
il est possible de d\'efinir un sous-groupe
$G^{0,-\infty}_{\mathcal{T}}(\bbS^{1})$ de $\Go(\bbS^{1})$ qui soit faiblement contractile.
Ce sous-groupe est obtenu en imposant des conditions suppl\'ementaires sur le symbole
principal, \`a savoir que ce dernier doit \^etre l'identit\'e sur 
\[
        \{s_{+}\}\sqcup \bbS^{1}_{-} \subset S^{*}\bbS^{1}
\] 
o\`u $s_{+}\in \bbS^{1}_{+}$ est un point de base choisi au pr\'ealable.  Avec ces restrictions,
$K^{0}(S^{*}\bbS^{1})$ est remplac\'e par
\[
    K^{0}(S^{*}\bbS^{1}, \{s_{+}\}\sqcup \bbS^{1}_{-})\cong \tilde{K}^{0}(\bbS^{1})\cong \{0\},
\]
alors que $K^{1}(S^{*}\bbS^{1})$ est remplac\'e par
\[
  K^{1}(S^{*}\bbS^{1}, \{s_{+}\}\sqcup \bbS^{1}_{-})\cong \tilde{K}^{1}(\bbS^{1})\cong \bbZ.
\]
Le th\'eor\`eme~\ref{gh.1} montre alors que les groupes d'homotopie de $G^{0,-\infty}_{\mathcal{T}}(\bbS^{1})$ sont tous triviaux.  En ce sens, le th\'eor\`eme~\ref{gh.1}
peut \^etre vu comme une g\'en\'eralisation du r\'esultat de contractibilit\'e de Melrose \cite{MelroseLMA}.

Une autre cons\'equence int\'eressante de la surjectivit\'e de l'homomorphisme de bord est la suivante.
\begin{theorem}
Soit $X$ une vari\'et\'e compacte sans bord de dimension sup\'erieure \`a z\'ero.
Soient $M$ un $CW$-complexe construit \`a partir d'un nombre fini
de cellules et $f:M\to \Gc(X)$ une application continue.  Si $i: \Gc(X)\hookrightarrow \Go(X)$
d\'enote l'inclusion canonique, alors l'application $i\circ f$ est homotope \`a l'application
constante
\[
    \begin{array}{lccc}
    \Id: & M & \to & \Go(X) \\
         & m & \mapsto & \Id
    \end{array}
\]
dans $\Go(X)$.
\label{gh.3}\end{theorem}
\begin{proof}
Puisque l'homomorphisme de bord $\pa: \pi_{k}(\Si)\to \pi_{k-1}(\Gc(X))$ est surjectif
pour tout $k\in \bbN$, on d\'eduit de la longue suite exacte de groupes d'homotopie que
\[
    i_{*}: \pi_{k}(\Gc(X))\to \pi_{k}(\Go(X))
\]
est une application triviale, c'est-\`a-dire envoie tout sur l'\'el\'ement identit\'e
de $\pi_{k}(\Go(X))$.  En utilisant la d\'ecomposition cellulaire de $M$, cela signifie que l'on
peut proc\'eder par r\'ecurrence pour construire une homotopie entre $i\circ f$ et 
$\Id: M\to \Go(X)$.
\end{proof}
\begin{remarque}
On dira que $\Gc(X)$ est \emph{faiblement r\'etractile} dans $\Go(X)$.
\label{gh.4}\end{remarque}

\section{Le cas des op\'erateurs suspendus}\label{os.0}

On peut obtenir un analogue des r\'esultats pr\'ec\'edents pour les op\'erateurs suspendus (suspended
operators en anglais) introduits par Melrose \cite{mr96h:58169}.  Rappelons d'abord bri\`evement leur
d\'efinition.  \`A nouveau, soit $X$ une vari\'et\'e lisse, compacte et sans bord et $l\in \bbN$ un
entier.  Consid\'erons l'espace $\Psi^{*}(X\times \bbR^{l})$ des op\'erateurs pseudodiff\'erentiels agissant sur la vari\'et\'e
non-compacte $X\times \bbR^{l}$.  Cet espace n'est pas une alg\`ebre, mais \`a tout le moins, chaque 
op\'erateur $A\in \Psi^{*}(X\times \bbR^{l})$ agit sur les fonctions lisses \`a support compact
\[
         A: \CI_{c}(X\times \bbR^{l}) \to \CI(X\times \bbR^{l}).
\]
Soit $T_{u}: X\times \bbR^{l}\to X\times \bbR^{l}$ le diff\'eomorphisme donn\'e par la translation 
$T_{u}(x,t)= (x,t+u)$ dans la deuxi\`eme variable et consid\'erons les op\'erateurs pseudodiff\'erentiels
dans $\Psi^{*}(X\times \bbR^{l})$ qui sont \textbf{invariants par translation}, c'est-\`a-dire satisfaisant
\begin{equation}
      T_{u}^{*}(Af)= A(T^{*}_{u}f), \quad \forall u\in \bbR^{l}, \; f\in \CI_{c}(X\times \bbR^{l}).
\label{os.1}\end{equation}
Le noyau de Schwartz $K_{A}$ d'un tel op\'erateur agit alors par convolution dans la direction de 
$\bbR^{l}$,
\[
              Af(x,t)= \int_{X}\int_{\bbR^{l}} K_{A}(x,x',t-s)f(x',s) ds
\]
avec $K_{A}\in \mathcal{C}^{-\infty}(X^{2}\times \bbR^{l}; \Omega_{R})$ o\`u $\Omega_{R}=\pi^{*}\Omega$ est
l'image r\'eciproque du fibr\'e des densit\'es sur $X$ par la projection 
\[
       \begin{array}{lccc}
       \pi: & X\times X\times \bbR^{l} & \to & X \\
            &            (x,x',t)   &\mapsto & x'.
       \end{array}
\]
Sous cette forme, ce noyau est alors singulier seulement sur la sous-vari\'et\'e 
$\{x=x', t=0\}$.  Pour pouvoir composer des op\'erateurs invariants par translation,
on peut aussi imposer une condition de d\'ecroissance rapide du noyau \`a l'infini
\begin{equation}
K_{A}\in \mathcal{C}_{c}^{-\infty}(X^{2}\times \bbR^{l}; \Omega_{R})+ \mathcal{S}(X^{2}\times \bbR^{l}; \Omega_{R}),
\label{os.2}\end{equation}
o\`u $\mathcal{S}(X^{2}\times\bbR^{l})$ d\'enote l'espace de Schwartz des sections \`a d\'ecroissance rapide
(avec toutes leurs d\'eriv\'ees) \`a l'infini.

\begin{definition}
Pour chaque $m\in \bbZ$ et $l\in \bbN$, on d\'efinit l'espace $\Psi^{m}_{\sus}(X)$
des op\'erateurs $l$ fois suspendus d'ordre $m$ 
sur $X$ comme \'etant le sous-espace de $\Psi^{m}(X\times \bbR^{l})$ constitu\'e des op\'erateurs
invariants par translation qui satisfont la condition de d\'ecroissance rapide 
\eqref{os.2}. 
\label{os.3}\end{definition}
Plus g\'en\'eralement, on peut d\'efinir l'espace des op\'erateurs suspendus agissant sur 
les sections d'un fibr\'e vectoriel complexe $E\to X$ par
\[
      \Psi^{m}_{\sus}(X;E):= \Psi^{m}_{\sus}(X)\otimes_{\CI(X^{2})} \CI(X^{2}; \Hom(E))
\]
o\`u $\Hom(E)$ est le fibr\'e vectoriel sur $X^{2}$ ayant pour fibre au-dessus de 
$(x,x')\in X^{2}$ 
\[
                 \Hom(E)_{(x,x')}= \hom(E_{x},E_{x'}).
\]
On peut v\'erifier qu'un op\'erateur suspendu $A\in \Psi^{m}_{\sus}(X;E)$ agit sur les
sections \`a d\'ecroissance rapide 
\[
        A: \mathcal{S}(X\times \bbR^{l};E) \to \mathcal{S}(X\times \bbR^{l};E)
\]
pour donner \`a nouveau des fonctions \`a d\'ecroissance rapide.  De l\`a, on peut voir que l'espace des 
op\'erateurs suspendus $\Psi^{*}_{\sus}(X;E)$ forme une alg\`ebre.  En prenant la 
transform\'ee de Fourier
\[
        K_{\hat{A}(\tau)}(y,y'):= \int_{\bbR^{l}} e^{-it\tau} K_{A}(y,y',t) dt, \quad \tau\in \bbR^{l},
\]
du noyau de Schwartz, on obtient une famille \`a $l$ param\`etres 
\[
     \hat{A}(\tau) \in \Psi^{m}(X;E), \tau \in \bbR^{l}
\]
d'op\'erateurs agissant sur $X$.  Cette famille est appel\'ee \textbf{famille indiciale} de $A$.  Un 
op\'erateur suspendu est compl\`etement d\'etermin\'e par cette derni\`ere et 
vice-versa.  En termes de la composition, un calcul direct montre que
\begin{equation}
      \widehat{A\circ B}(\tau) = \hat{A}(\tau)\circ \hat{B}(\tau), \quad \forall \tau\in \bbR^{l}.
\label{os.4}\end{equation}
Par cons\'equent, on voit que pour un op\'erateur suspendu inversible $A$, sa famille
indiciale $\hat{A}(\tau)$ est aussi inversible pour tout $\tau$.  La r\'eciproque est aussi vraie,
mais moins \'evidente (voir \cite{mr96h:58169}).  En quelque sorte, la famille indiciale 
peut \^etre vue comme un symbole dans la variable $\tau\in \bbR^{l}$ qui est
toujours quantifi\'e dans les variables $(x,\xi)\in T^{*}X$.  

Les op\'erateurs suspendus ont aussi un symbole principal qui donne lieu \`a une suite exacte
\begin{equation}
0\to \Psi_{\sus}^{m-1}(X;E)\to \Psi^{m}_{\sus}(X;E) \overset{\sigma_{m}}{\to}
\CI(S^{*}_{X}(X\times \bbR); \pi^{*}\hom E\otimes D_{m})\to 0  
\label{os.5}\end{equation}
o\`u $S^{*}(X\times \bbR)= T^{*}(X\times \bbR)\setminus0 /\bbR^{+}$ et $S^{*}_{X}(X\times \bbR)$ est
sa restriction \`a $X\times\{ 0\}$, alors que $D_{m}$ est le fibr\'e en droite sur 
$S^{*}(X\times \bbR)$ d\'efini par les fonctions homog\`enes de degr\'e $m$ sur $T^{*}(X\times\bbR)
\setminus 0$.
\begin{definition}
Un op\'erateur suspendu $A\in \Psi^{m}(X;E)$ est \emph{elliptique} si son symbole principal
est inversible. 
\label{os.6}\end{definition}
On peut v\'erifier via la construction d'un inverse modulo $\Psi^{-\infty}_{\sus}(X;E)$ 
que la famille indiciale $\hat{A}(\tau)$ d'un op\'erateur suspendu elliptique 
$A\in \Psi_{\sus}^{m}(X;E)$ est inversible pour tout $\tau\in \bbR^{l}$ tel que $|\tau|>R$, o\`u 
$R>0$ est choisi suffisamment grand.

On peut maintenant d\'efinir la version stabilis\'ee du groupes des op\'erateurs suspendus 
inversibles d'ordre z\'ero par 
\begin{equation}
   \Gos(X):= \{\Id+Q\quad | \quad Q\in \CI(\bbS^{1}\times \bbS^{1}; \Psi^{0}_{\sus}(X)),
     \; \Id+Q \;\mbox{est inversible} \}
\label{os.7}\end{equation}
o\`u $Q\in \CI(\bbS^{1}\times \bbS^{1}; \Psi^{0}_{\sus}(X))$ agit sur $f\in 
\mathcal{S}(X\times \bbR^{l}\times \bbS^{1})$ par
\[
             (Qf)(x,t,\theta)= \int_{\bbS^{1}} (Q(\theta,\theta')f_{\theta'})(x,t) d\theta',
\]
la fonction $f_{\theta'}\in \mathcal{S}(X\times \bbR^{l})$ \'etant d\'efinie par $f_{\theta'}(x,t):= 
f(x,t,\theta)$.  De m\^eme, on peut d\'efinir le sous-groupe des 
\flqq perturbations compactes inversibles\frqq \ de l'identit\'e par 
\begin{equation}
G^{-1}_{\sus}(X):= \{  \Id+Q\in G^{0}_{\sus}(X)\quad | \quad 
    Q\in \CI(\bbS^{1}\times \bbS^{1}; \Psi^{-1}_{\sus}(X))\},
\label{os.8}\end{equation} 
en ce sens que la famille indiciale d'un op\'erateur $A\in \Psi_{\sus}^{-1}(X)$ est
constitu\'ee d'op\'erateurs compacts. 
\begin{lemme}
Le symbole principal donne lieu \`a une suite exacte
\[
0\to \Gcs(X)\to\Gos(X)\overset{\sigma}{\to} \Ss\to 0 
\]
avec 
\[
\Ss:= \left\{ \begin{array}{cl}
            \CI(S^{*}_{X}(X\times\bbR^{l});G^{-\infty}(\bbS^{1})), & l \; \mbox{impair},  \\
            \ker[\widehat{\ind}:\CI(S^{*}_{X}(X\times\bbR^{l});G^{-\infty}(\bbS^{1})) \to \bbZ], &
                                                                     l \; \mbox{pair},
              \end{array} \right.
\]
 o\`u $\widehat{\ind}: \CI(S^{*}_{X}(X\times\bbR^{l});G^{-\infty}(\bbS^{1}))\to \bbZ$ est un indice de famille 
 d\'efini \`a partir de la famille indiciale (voir \eqref{ind.1} plus bas).
\label{os.9}\end{lemme}
\begin{proof}
\'Etant donn\'e un 
symbole $a\in \Ss$, on peut trouver un op\'erateur 
$A\in \CI(\bbS^{1}\times \bbS^{1}; \Psi^{0}_{\sus}(X))$ tel que
\[
     \sigma(A)=a.
\]
Cet op\'erateur $A$ est donc elliptique et par cons\'equent sa famille indiciale
$\hat{A}(\tau)$ est inversible pour tout $\tau$ satisfaisant $|\tau|>R$ o\`u $R$ est 
une constante positive assez grande.  Vue comme une famille d'op\'erateurs de Fredholm,
c'est dire que la famille indiciale $\hat{A}$ d\'efinit un indice de famille 
\begin{equation}
    \ind\hat{A}\in K^{0}_{c}(\bbR^{l})\cong \tilde{K}^{0}(\bbS^{l})\cong 
    \left\{ \begin{array}{cl}
              \{0\}, & l \; \mbox{impair},  \\
              \bbZ,  & l \; \mbox{pair}.
            \end{array} \right.
\label{ind.2}\end{equation}
Ainsi, lorsque $l$ est impair, cet indice est n\'ecessairement trivial et il n'y a aucune
obstruction \`a l'existence d'un op\'erateur $Q\in \CI(\bbS^{1}\times \bbS^{1}, 
\Psi^{-\infty}_{\sus}(X))$ tel que $\hat{A}(\tau)+ \hat{Q}(\tau)$ soit inversible
pour tout $\tau\in \bbR^{l}$.  On en d\'eduit que $A+Q\in \Gos(X)$ avec $\sigma(A+Q)=a$,
d'o\`u la surjectivit\'e \`a droite de la suite exacte lorsque $l$ est impair.

Lorsque $l$ est pair, l'indice \eqref{ind.2} ne d\'epend pas du choix de $A$ tel
que $\sigma(A)=a$ et donc d\'efinit une application
\begin{equation}
\widehat{\ind}: \CI(S^{*}_{X}(X\times \bbR^{l});G^{-\infty}(\bbS^{1}))\to \bbZ.
\label{ind.1}\end{equation}
Comme dans le lemme~\ref{fdS.4}, cet indice  mesure exactement l'obstruction
\`a trouver un op\'erateur $A\in \Gos(X)$ tel que $\sigma(A)=a$, ce qui \'etablit 
la surjectivit\'e \`a droite de la suite exacte dans le cas o\`u $l$ est pair en
posant
\[
  \Ss:= \ker[\widehat{\ind}: K^{1}(S^{*}_{X}(X\times \bbR^{l}))\to \bbZ].
\]

\end{proof}

\`A nouveau, on peut v\'erifier que la suite exacte du lemme pr\'ec\'edent est une 
fibration de Serre, ce qui donne une longue suite exacte de groupes d'homotopie
\begin{multline}
\cdots\to \pi_{k}(\Gcs(X))\to \pi_{k}(\Gos(X))\to \pi_{k}(\Ss)\overset{\pa}{\to}
  \pi_{k-1}(\Gcs(X))\to \cdots \\
\cdots\to \pi_{1}(\Ss)\overset{\pa}{\to} \pi_{0}(\Gcs(X))\to \pi_{0}(\Gos(X))\to
  \pi_{0}(\Ss).
\label{os.10}\end{multline}
Pour $p\in \bbN$, $\Gc_{\susn{2p}}(X)$ est un espace 
classifiant pour la $K$-th\'eorie \textbf{impaire}, alors que $\Gc_{\susn{2p-1}}(X)$ est
un espace classifiant pour la $K$-th\'eorie \textbf{paire}. 
Utilisant le fait que $G^{-\infty}(\bbS^{1})$ est un espace classifiant pour la 
$K$-th\'eorie impaire, on a aussi que pour $k>0$,
\begin{equation}
\begin{aligned}
\pi_{k}(\Ss)&\cong [ S^{k}(S^{*}_{X}(X\times \bbR^{l})); G^{-\infty}(\bbS^{1})] \\
            &\cong K^{-k-1}(S^{*}_{X}(X\times \bbR^{l})) \\
            &\cong \left\{
                   \begin{array}{ll}
                   K^{0}(S^{*}_{X}(X\times \bbR^{l})), & k \, \mbox{impair}, \\
                   K^{1}(S^{*}_{X}(X\times \bbR^{l}), & k \, pair.
                   \end{array}\right. 
\end{aligned}
\label{os.11}\end{equation}
Ce r\'esultat est aussi valable pour $k=0$ et $l$ impair, mais pour $k=0$ et $l$ pair, on
a plut\^ot
\[
   \pi_{0}(\Ss)\cong \ker[\widehat{\ind}:K^{1}(S^{*}_{X}(X\times \bbR^{l}))\to \bbZ]                   
\]
l'indice \eqref{ind.1} ne d\'ependant que de la $K$-classe d\'efinie par le symbole principal.

On pourra donc calculer les groupes d'homotopie de $\Gos(X)$ en montrant que l'homomorphisme
de bord $\pa$ est surjectif.  La preuve est tr\`es similaire au cas des op\'erateurs pseudodiff\'erentiels usuels.  Dans un premier temps, on montre que l'homomorphisme de bord correspond
\`a un indice de famille.  On montre alors que cet indice de famille est surjectif en utilisant
le th\'eor\`eme d'Atiyah-Singer \cite{Atiyah-Singer4} pour les familles d'op\'erateurs.  

L'indice de famille qu'il faut consid\'erer est obtenu en regardant un op\'erateur
suspendu elliptique comme une famille \`a $l$ param\`etres d'op\'erateurs de Fredholm inversibles
\`a l'infini.  Plus pr\'ecis\'ement, soit $f:\bbS^{k}\to \Ss$ une application
repr\'esentant un \'el\'ement de $\pi_{k}(\Ss)$.
Sans perte de g\'en\'eralit\'e, on peut supposer que $f\equiv \Id$ dans un voisinage du point de 
base $s_{0}$ de $\bbS^{k}$.  
Soit alors 
\[
    \tilde{f}: \bbS^{k} \to \Id + \CI(\bbS^{1}\times \bbS^{1}, \Psi^{0}_{\sus}(X))
\]
un rel\`evement de $f$ dans $\Id+ \CI(\bbS^{1}\times \bbS^{1}, \Psi^{0}_{\sus}(X))$ tel
que $\tilde{f}\equiv \Id$ dans un voisinage du point de base de $\bbS^{k}$.  Alors la famille
indiciale de $\tf$, d\'enot\'ee $\hf$, d\'efinit une famille d'op\'erateurs de Fredholm
\[
     \hf: \bbS^{k}\times \bbR^{l} \to \cF(\cH), \quad \cH:= \Ld(X\times \bbS^{1}),
\]
inversible \`a l'infini et sur $\{s_{0}\}\times \bbR^{l}$ o\`u $s_{0}\in \bbS^{k}$ est le point 
de base de $\bbS^{k}$.  Cela d\'efinit donc un indice de famille
\[
   \ind(\hf)\in \tilde{K}^{0}(\bbS^{k+l}).
\]
En effet, la condition d'inversibilit\'e assure que $\hf$ d\'efinit un indice sur la 
$l$-suspension de $\bbS^{k}$ en identifiant $\bbR^{l}\cup \{\infty\}$ avec $\bbS^{l}$, 
$\{\infty\}$ \'etant le point de base.
Cet indice ne d\'epend pas du choix du rel\`evement $\hf$ ou du choix
$f$ du repr\'esentant de la classe d'homotopie.  On a donc en fait une application 
\begin{equation}
    \ind: \pi_{k}(\Ss)\to \tilde{K}^{0}(\bbS^{k+l}).
\label{os.14}\end{equation}
D'autre part, on a la s\'erie d'identifications
\begin{equation}
\begin{aligned}
   \pi_{k-1}(\Gcs(X)) & \cong \pi_{k+l-1}(\Gc(X))\cong \tilde{K}^{-1}(\bbS^{k+l-1})  \\
                      & \cong \tilde{K}^{0}(\bbS^{k+l}).
\end{aligned}
\label{os.12}\end{equation}
\begin{lemme}
Soit $p:\pi_{k-1}(\Gcs(X))\to\tilde{K}^{0}(\bbS^{k+l})$ l'isomorphisme r\'esultant
de l'identification \eqref{os.12}, alors l'homomorphisme de bord est donn\'e par
\[
    \pa= p^{-1}\circ \ind
\]
o\`u $\ind$ est l'indice de famille \eqref{os.14}.
\label{os.13}\end{lemme}
\begin{proof}
Modulo quelques adaptations mineures, la d\'emonstration est la m\^eme que celle de la 
proposition~\ref{chb.3}.  On laisse le soin au lecteur de compl\'eter les d\'etails. 
\end{proof}
\begin{lemme}
L'homomorphisme de bord $\pa: \pi_{k}(\Ss)\to \pi_{k-1}(\Gcs(X))$ est surjectif.
\label{os.15}\end{lemme}
\begin{proof}
Par le lemme pr\'ec\'edent, il suffit de montrer que l'indice de famille 
\[
   \ind: \pi_{k}(\Ss) \to \tilde{K}^{0}(\bbS^{k+l})
\]
est surjectif.  Lorsque $k+l$ est impair, $\tilde{K}^{0}(\bbS^{k+l})\cong \{0\}$ et 
il n'y a rien \`a montrer.  Lorsque $k+l$ est pair, on proc\`ede comme dans le cas
$k=0$ du lemme~\ref{chb.4}, mais cette fois en utilisant l'indice 
d'Atiyah-Singer~\cite{Atiyah-Singer4} pour les familles d'op\'erateurs.
\end{proof}

\begin{theorem}
Soit $X$ une vari\'et\'e compacte sans bord de dimension sup\'erieure \`a z\'ero.
Lorsque $l$ est \textbf{impair}, les groupes d'homotopie de $\Gos(X)$ sont donn\'es par 
\[
          \pi_{k}(\Gos(X))\cong \left\{
             \begin{array}{ll}
             K^{-1}(S^{*}_{X}(X\times\bbR^{l})), & k\, \mbox{pair},  \\ 
               \ker[ \ind_{t}\circ \delta_{l}; K^{0}(S^{*}_{X}(X\times \bbR^{l}))\to \bbZ], &
                                 k\, \mbox{impair},                                
             \end{array}\right.
\]
o\`u $\delta_{l}: K^{0}(S^{*}_{X}(X\times\bbR^{l}))\to K_{c}^{-1}(T^{*}X\times \bbR^{l})\cong
K^{0}(T^{*}X)$ est l'homomorphisme de bord (compos\'e avec la p\'eriodicit\'e de Bott) associ\'e
\`a la paire $(\overline{T^{*}X\times\bbR^{l}}, S^{*}_{X}(X\times\bbR^{l}))$, alors que $ind_{t}$ est 
l'indice topologique de Atiyah-Singer.

Lorsque $l$ est \textbf{pair}, les groupes d'homotopie de $\Gos(X)$ sont plut\^ot donn\'es par
\[
          \pi_{k}(\Gos(X))\cong \left\{
             \begin{array}{ll}
             \ker[\ind_{t}\circ\delta_{l}:K^{-1}(S^{*}_{X}(X\times\bbR^{l}))\to \bbZ], & k\, \mbox{pair},  \\ 
               K^{0}(S^{*}_{X}(X\times \bbR^{l})), &
                                 k\, \mbox{impair},                                
             \end{array}\right.
\]
o\`u $\delta_{l}: K^{-1}(S^{*}_{X}(X\times \bbR^{l}))\to K^{0}_{c}(T^{*}X\times \bbR^{l})\cong 
K^{0}_{c}(T^{*}X)$ est l'homomorphisme de bord associ\'e \`a la paire 
$(\overline{T^{*}X\times\bbR^{l}}, S^{*}_{X}(X\times\bbR^{l}))$.

\label{os.16}\end{theorem}
\begin{proof}
C'est une cons\'equence de la surjectivit\'e de l'homomorphisme de bord et de la longue suite
exacte~\eqref{os.10}.  Dans le cas $k=0$, le r\'esultat est une cons\'equence de la surjectivit\'e
\`a droite de la suite exacte du lemme~\ref{os.9}.
\end{proof}

Comme dans le cas des op\'erateurs diff\'erentiels usuels, la surjectivit\'e de l'homomorphisme
de bord a aussi la cons\'equence suivante.

\begin{theorem}
Soit $X$ une vari\'et\'e compacte sans bord de dimension sup\'erieure \`a z\'ero.  Alors
le sous-espace $\Gcs(X)$ est faiblement r\'etractile dans $\Gos(X)$ (voir la remarque~\ref{gh.4}). 
\label{os.17}\end{theorem}
\begin{proof}
La d\'emonstration est la m\^eme que dans le th\'eor\`eme~\ref{gh.3}.
\end{proof}

\section{Une Application reli\'ee \`a la th\'eorie de l'indice} \label{ati.0} 

Ce dernier r\'esultat donne lieu \`a une application int\'eressante en th\'eorie de l'indice pour les op\'erateurs \`a cusp fibr\'e (fibred cusp operators
dans la terminologie de Melrose).  Ceux-ci ont \'et\'e introduits par
Mazzeo et Melrose dans \cite{mazzeo-melrose4}.  Par courtoisie pour le 
lecteur, rappelons d'abord bri\`evement leur d\'efinition et leurs 
propri\'et\'es de base.  On r\'ef\`ere \`a l'article de Mazzeo et Melrose
\cite{mazzeo-melrose4} ainsi qu'\`a \cite{Lauter-Moroianu} et 
\cite{bpfco} pour de plus amples d\'etails.

Soit $M$ une vari\'et\'e compacte avec bord $\pa M$ de classe $\CI$.  Supposons que son bord 
soit muni d'une structure de fibration de classe $\CI$ localement triviale,
\begin{equation}
\xymatrix{ Z \ar@{-}[r]  &  \pa M\ar[d]^{\Phi} \\
  &  Y }
\label{ti.1}\end{equation}
et soit $x\in \CI(M)$ une fonction de d\'efinition du bord, c'est-\`a-dire 
que $\pa M= \{x=0\}$, $x>0$ sur $M\setminus \pa M$ et la diff\'erentielle
$dx$ est diff\'erente de z\'ero partout sur le bord $\pa M$.  
\begin{definition}
Un \textbf{champ de vecteur \`a cusp fibr\'e} $V\in \CI(M;TM)$ est un
 champ de vecteur sur $M$ qui est tangent aux fibres de la fibration 
$\Phi:\pa M\to Y$ et tel que $Vx\in x^{2}\CI(M)$.  On d\'enote par
$\cV_{\Phi}(M)$ l'ensemble des champs de vecteurs \`a cusp fibr\'e de $M$.  
\label{fc.1}\end{definition}
Soient $(x,y,z)$ des coordonn\'ees dans un voisinage de $p\in \pa M$, 
$x$ \'etant la fonction de d\'efinition du bord et $y$ et $z$ \'etant des
coordonn\'ees locales sur $Y$ et $Z$ respectivement, o\`u on assume que la
fibration $\Phi$ est triviale dans le voisinage du point $p$ consid\'er\'e.  
Dans ces coordonn\'ees, un champ de vecteur a cusp fibr\'e prend la forme 
suivante
\begin{equation}
   V= a x^{2}\frac{\pa}{\pa x} + \sum_{i=1}^{l} b_{i} x\frac{\pa}{\pa y_{i}}
     + \sum_{i=1}^{k}c_{i} \frac{\pa}{\pa z_{i}}
\label{fc.2}\end{equation}
avec $l=\dim Y$, $k=\dim Z$ et o\`u $a, b_{i}, c_{i}$ sont des fonctions de 
classe $\CI$.  

En utilisant le crochet de Lie pour les champs de vecteur, un calcul
direct montre que $\cV_{\Phi}(M)$ est en fait une alg\`ebre de Lie.  On 
d\'efinit l'espace des op\'erateurs diff\'erentiels \`a cusp fibr\'e comme
\'etant la $\CI(M)$-alg\`ebre universelle enveloppante de l'alg\`ebre de 
Lie $\cV_{\Phi}(M)$.  L'espace des op\'erateurs diff\'erentiels \`a cusp fibr\'e d'ordre
$m$ est donc l'espace des op\'erateurs diff\'erentiels engendr\'es par 
$\CI(M)$ et la composition d'au plus $m$ \'el\'ements de $\cV_{\Phi}(M)$.
Dans les coordonn\'ees locales $(x,y,z)$, un op\'erateur \`a cusp fibr\'e
$P$ d'ordre $m$ est de la forme
\begin{equation}
       P= \sum_{|\alpha| + |\beta| + q \le m} p_{\alpha,\beta,q} 
              \left( x^{2}\frac{\pa}{\pa x}\right)^{q}\left( x\frac{\pa}{\pa y}
                 \right)^{\alpha} \left( \frac{\pa}{\pa z}\right)^{\beta}
\label{fc.3}\end{equation}
o\`u les $p_{\alpha,\beta, q}$ sont des fonctions de classe $\CI$.

Plus g\'en\'eralement et \`a un niveau intuitif, on obtient un op\'erateur
\textbf{pseudodiff\'erentiel} \`a cusp fibr\'e $P$ en permettant \`a $P$
d'\^etre une \flqq fonction\frqq  plus g\'en\'erale qu'un polynome dans les 
 \flqq variables\frqq  $x^{2}\frac{\pa}{\pa x}, x\frac{\pa}{\pa y}$ et  
$\frac{\pa}{\pa z}$.  Sans entrer dans les d\'etails, disons qu'une 
d\'efinition rigoureuse de ces op\'erateurs passe par la description 
de leurs noyaux de Schwartz, qui sont des distributions sur une version
\'eclat\'ee (au sens de Melrose) du produit cart\'esien $M\times M$.
On d\'enote par $\Psi^{m}_{\Phi}(M;E,F)$ l'espace des op\'erateurs 
pseudodiff\'erentiels \`a cusp fibr\'e d'ordre $m$ agissant sur les sections
du fibr\'e vectoriel $E\to M$ pour donner des sections du fibr\'e vectoriel
$F\to M$.  

\'Etant donn\'es des fibr\'es vectoriels complexes $E$, $F$ et $G$ sur $M$, les op\'erateurs
\`a cusp fibr\'e satisfont la propri\'et\'e de composition usuelle
\[
    \Psi_{\Phi}^{m}(M;F,G) \circ \Psi_{\Phi}^{n}(M;E,F) \subset \Psi_{\Phi}^{m+n}(M;E,G).
\]
De plus, lorsqu'un op\'erateur \`a cusp fibr\'e est inversible, son inverse est aussi un op\'erateur
\`a cusp fibr\'e.  Dans \cite{mazzeo-melrose4}, Mazzeo et Melrose donne une caract\'erisation \'el\'egante et pr\'ecise des op\'erateurs \`a cusp fibr\'e qui sont de Fredholm lorsqu'agissant sur des espaces
de Sobolev adapt\'es \`a ce contexte.  Il y a d'abord une notion de symbole principal.  Celle-ci est 
obtenue en rempla\c{c}ant le fibr\'e tangent $TM$ par le fibr\'e tangent \`a cusp fibr\'e 
${}^{\Phi}TM\to M$, qui est d\'efini de sorte qu'on ait une identification canonique
\[
             \CI(M; {}^{\Phi}TM)= \cV_{\Phi}(M).
\]
Le fibr\'e ${}^{\Phi}TM$ est isomorphe au fibr\'e tangent $TM$, mais pas d'une fa\c{c}on canonique.  Soit
alors ${}^{\Phi}S^{*}M:= ({}^{\Phi}TM\setminus M)/ \bbR^{+}$ le fibr\'e cosph\'erique \`a cusp fibr\'e
de $M$, o\`u ${}^{\Phi}T^{*}M$ est le fibr\'e dual de ${}^{\Phi}TM$.  Soit aussi $\mathcal{R}^{m}\to {}^{\Phi}S^{*}M$ le fibr\'e trivial en droites complexes dont les sections sont donn\'ees par les fonctions
sur ${}^{\Phi}T^{*}M\setminus M$ qui sont homog\`enes de degr\'e $m$.  Il est alors possible de d\'efinir
une application symbolique
\[
      \sigma_{m}: \Psi_{\Phi}^{m}(M;E,F) \to \CI({}^{\Phi}S^{*}M; \hom(E,F)\otimes \mathcal{R}^{m})
\]
qui donne lieu \`a la suite exacte
\[
  0\to \Psi^{m-1}_{\Phi}(M;E,F)\longrightarrow \Psi_{\Phi}^{m}(M;E,F) 
  \overset{\sigma_{m}}{\longrightarrow}\CI({}^{\Phi}S^{*}M;
  \hom(E,F)\otimes \mathcal{R}^{m})\to 0. 
\]
\begin{definition}
On dit qu'un op\'erateur $P\in \Psi^{m}_{\Phi}(M;E,F)$ est \textbf{elliptique} si son symbole 
principal $\sigma_{m}(P)$ est inversible.
\label{fc.4}\end{definition}
Vu la pr\'esence d'un bord $\pa M$, il n'est pas suffisant de demander qu'un op\'erateur soit
elliptique pour qu'il soit du m\^eme coup de Fredholm.  Cela d\'ecoule du fait qu'un op\'erateur
$P\in \Psi^{-1}_{\Phi}(M;E,F)$ d'ordre $-1$ n'est pas n\'ecessairement compact.  Cet op\'erateur 
doit aussi satisfaire une certaine condition de d\'ecroissance au bord.  Mazzeo et Melrose ont 
d\'ecrit cette condition gr\^ace \`a l'introduction d'une autre suite exacte,
\[
0\to x\Psi_{\Phi}^{m}(M;E,F)\longrightarrow \Psi_{\Phi}^{m}(M;E,F)
\overset{N_{\Phi}}{\longrightarrow} \Psi_{\Phi-\susf}^{m}(\pa M;E,F)\to 0,
\]
o\`u $N_{\Phi}:\Psi^{m}_{\Phi}(M;E,F)\to \Psi^{m}_{\Phi-\susf}(\pa M;E,F)$ est \textbf{l'application normale}
et $N_{\Phi}(P)$ est appel\'e \textbf{l'op\'erateur normal} de $P$ (parfois aussi la famille indiciale de $P$).  Soit ${}^{\Phi}NY\cong TY\times \bbR$ le fibr\'e noyau associ\'e \`a la restriction canonique
\[
{}^{\Phi}T_{\pa M}M\to T_{\pa M}M.
\]
Alors $\Psi^{m}_{\Phi-\susf}(\pa M;E,F)$ 
d\'enote l'espace des op\'erateurs ${}^{\Phi}NY$-suspendus d'ordre $m$.   Un \'el\'ement de 
$\Psi^{m}_{\Phi-\susf}(\pa M;E,F)$ est une famille d'op\'erateur ($l+1$)-suspendus param\'etris\'ee par
$Y$ o\`u $l:= \dim Y$.  Pour un $y\in Y$ fix\'e, un \'el\'ement de $\Psi^{m}_{\Phi-\susf}(\pa M;E,F)$ 
d\'efinit donc un op\'erateur dans $\Psi_{\susn{l+1}}^{m}(\Phi^{-1}(y);E,F)$ lorsqu'on identifie 
la fibre ${}^{\Phi}N_{y}Y$ avec $\bbR^{l+1}$.  L'op\'erateur normal joue en quelque sorte le r\^ole d'un
symbole principal qui caract\'erise le comportement asymptotique de l'op\'erateur pr\`es du bord.  
L'op\'erateur normal est symbolique et local par rapport \`a la variable $y\in Y$, mais est toujours
quantifi\'e par rapport \`a la variable $z\in Z$.  L'application normale est un homomorphisme d'alg\`ebre,
\[
    N_{\Phi}(A\circ B)= N_{\Phi}(A)\circ N_{\Phi}(B)
\]
pour $A\in \Psi^{k}_{\Phi}(M;F,G)$ et $B\in \Psi_{\Phi}^{m}(M;E,F)$.  
\begin{definition}
Un op\'erateur \`a cusp fibr\'e $A\in \Psi_{\Phi}^{m}(M;E,F)$ est dit \textbf{totalement elliptique}
(fully elliptic dans la terminologie de Melrose)
s'il est elliptique et si son op\'erateur normal $N_{\Phi}(A)$ est inversible.
\label{fc.5}\end{definition}
Soit $\Ld(M;E)$ l'espace des sections du fibr\'e $E\to M$ de carr\'e absolument int\'egrable par 
rapport \`a un choix de m\'etrique sur $E$ et un choix de forme volume sur $M$.  Pour $l\in \bbR$ et
$m\in \bbN_{0}$, consid\'erons les espaces de Sobolev
\begin{equation}
\begin{gathered}
x^{l} H^{m}_{\Phi}(M;E):= \{ u\in x^{l}\Ld(M;E) \quad | \quad Pu \in x^{l}\Ld(M;E)\; \forall P\in \Psi^{m}_{\Phi}(M;E)\}, \\
x^{l} H^{-m}_{\Phi}(M;E):= \{ u \in \mathcal{C}^{-\infty}(M;E) \quad | \quad 
                   u\in \sum_{i=1}^{N} P_{i}u_{i}, \hspace{3.5cm}\\
                   \hspace{5cm} u_{i}\in x^{l}\Ld(M;E), P_{i}\in \Psi^{m}_{\Phi}(M;E) \}.
\end{gathered}
\label{fc.6}\end{equation} 
Pour ces espaces de Sobolev, Mazzeo et Melrose ont donn\'e la caract\'erisation suivante des op\'erateurs
\`a cusp fibr\'e qui sont de Fredholm.
\begin{proposition}
Pour $m\in \bbZ$, un op\'erateur $P\in \Psi^{m}_{\Phi}(M;E,F)$ donne lieu \`a un op\'erateur
born\'e
\[
         P: x^{l}H^{m'}_{\Phi}(M;E)\to x^{l}H^{m'-m}_{\Phi}(M;F)
\]
qui est de Fredholm si et seulement si $P$ est totalement elliptique.
\label{fc.7}\end{proposition}
L'\'etude de l'indice des op\'erateurs \`a cusp fibr\'e qui sont de Fredholm a \'et\'e l'objet
de plusieurs travaux depuis la parution de l'article de Mazzeo et Melrose \cite{mazzeo-melrose4},
notamment par Nye et Singer \cite{Nye-Singer}, Lauter et Moroianu \cite{Lauter-Moroianu}, \cite{Lauter-Moroianu2}, \cite{Lauter-Moroianu3}, Melrose et l'auteur \cite{iKfaficu}, Moroianu \cite{Moroianu2} ainsi que Leichtnam, Mazzeo et Piazza \cite{Leichtnam-Mazzeo-Piazza}.

Rappelons que la d\'efinition des op\'erateurs \`a cusp fibr\'e donne lieu \`a deux extr\^emes.  D'abord,
on a le cas o\`u la base $Y$ de la fibration $\Phi: \pa M\to Y$ n'est constitu\'ee que d'un seul point.
Les op\'erateurs associ\'es \`a une telle fibration sont appel\'es op\'erateurs cuspidaux (cusp
operators dans la terminologie de Melrose).  L'op\'erateur normal $N_{\Phi}(P)$ d'un op\'erateur cuspidal
$P\in \Psi_{\Phi}^{m}(M;E,F)$ est un \'el\'ement de $\Psi_{\susn{1}}^{m}(M;E,F)$, c'est-\`a-dire que c'est 
un op\'erateur suspendu dans le sens usuel.  Intuitivement, $N_{\Phi}(P)$ peut \^etre interpr\'et\'e comme la partie de l'op\'erateur $P$ asymptotiquement invariante par translation dans la direction normale au bord.  Dans \cite{Melrose-Nistor}, Melrose et Nistor ont obtenu
une formule pour l'indice des op\'erateurs cuspidaux qui g\'en\'eralise le th\'eor\`eme d'indice de 
Atiyah, Patodi et Singer \cite{APS1}.  L'autre cas extr\^eme est obtenu en demandant que $Y=\pa M$ et 
que la fibration soit simplement donn\'ee par l'identit\'e $\Id: \pa M\to \pa M$.  Les op\'erateurs 
qu'on obtient ainsi sont appel\'es op\'erateurs de diffusion (scattering operators dans la terminologie
de Melrose).  Comme la fibre de la fibration est de dimension nulle, l'op\'erateur normal est 
compl\`etement symbolique pour un op\'erateur de diffusion.
Dans ce cas, l'indice peut-\^etre d\'ecrit \`a l'aide d'un \'el\'ement de $K^{0}_{c}(T^{*}M,
\left. T^{*}M\right|_{\pa M})$ et la formule qu'on obtient (voir \cite{MelroseGST}, th\'eor\`eme 6.4), est tr\`es similaire \`a celle du th\'eor\`eme
d'Atiyah et Singer.  

Le cas g\'en\'eral est \`a mi-chemin entre ces deux extr\^emes.  Lorsqu'on veut \'etudier l'indice des
op\'erateurs \`a cusp fibr\'e, une des difficult\'es importantes est de comprendre l'int\'eraction subtile
qu'il y a entre le symbole principal et l'op\'erateur normal.  On en vient au r\'esultat principal de ce 
paragraphe.  Il consiste \`a montrer que, \'etant donn\'e un op\'erateur totalement elliptique
$P\in \Psi^{m}_{\Phi}(M;E,F)$, on peut toujours choisir modulo certaines d\'eformations comment quantifier
le symbole de l'op\'erateur normal de $P$ en demandant que celui-ci soit inversible et sans changer l'indice
de $P$.

\begin{corollaire}
Soit $P\in \Psi_{\Phi}^{0}(M;E,F)$ un op\'erateur \`a cusp fibr\'e totalement elliptique.  Si un \'el\'ement inversible $p_{1}\in \Psi^{0}_{\Phi-\susf}(\pa M;E,F)$ est tel que son symbole 
principal $\sigma_{0}(p_{1})$ est \'egal \`a celui de 
\[
p_{0}:=N(P)\in \Psi_{\Phi-\susf}(\pa M;E,F),
\]
l'op\'erateur normal de $P$, alors, possiblement apr\`es stabilisation des fibr\'es vectoriels
$E$ et $F$, il existe une famille lisse \`a un param\`etre $P_{t}\in \Psi^{0}_{\Phi}(M;E,F)$,
$t\in [0,1]$, d'op\'erateurs totalement elliptiques tels que $P_{0}=P$ et $N(P_{1})=p_{1}$.
\label{os.18}\end{corollaire}
\begin{proof}
Remarquons d'abord que $p_{0}^{-1}\circ p_{1}\in G^{-1}_{\Phi-\susf}(\pa M;E)$.  Or, dans un 
ouvert $\mathcal{U}$ de $Y$ o\`u la fibration $\Phi$ et le fibr\'e vectoriel $T^{*}\mathcal{U}\to
\mathcal{U}$ sont triviaux, un op\'erateur 
$A\in \Psi^{0}_{\Phi-\susf}( \Phi^{-1}(\mathcal{U});E)$ peut \^etre vu comme une famille d'op\'erateurs
$(l+1)$-suspendus 
\[
     A: \mathcal{U}\to \Psi^{0}_{\susn{l+1}}(Z;E)
\]
o\`u $l=\dim Y$.  En choisissant une d\'ecomposition cellulaire de $Y$ telle que chaque
cellule soit contenue dans un ouvert trivialisant \`a la fois la fibration $\Phi$ et 
le fibr\'e $T^{*}{Y}$, on peut alors proc\'eder par r\'ecurrence en utilisant le th\'eor\`eme
pr\'ec\'edent pour construire une homotopie entre $p_{0}^{-1}p_{1}$ et l'identit\'e dans
$G^{0}_{\Phi-\susf}(\pa M;E)$ (possiblement en stabilisant $E$).  Par composition avec $p_{0}$, cela donne une homotopie entre $p_{0}$ et $p_{1}$ donn\'ee par
$p_{t}\in \Psi^{0}_{\Phi-\susf}(\pa M;E,F)$ inversible pour tout $t\in [0,1]$.  Il est alors facile
de relever cette homotopie parmi les op\'erateurs totalment elliptiques pour obtenir le r\'esultat 
d\'esir\'e.   
\end{proof}

Ce r\'esultat permet de r\'eduire l'\'etude de l'indice \`a des op\'erateurs totalement elliptiques
$P$ dont l'op\'erateur normal sastifait certaines conditions, par exemple tel que le symbole total
\footnote{Bien d\'efini apr\`es avoir fix\'e un choix de quantification.} soit identiquement nul
sauf pour le symbole principal $\sigma(N_{\Phi}(P))$.

\section{Une description topologique du d\'eterminant r\'esiduel} \label{dr.0}

Le d\'eterminant r\'esiduel (residue determinant en anglais) a \'et\'e introduit par Simon Scott dans \cite{Scott}.  C'est une fonctionnelle qui joue le r\^ole de d\'eterminant pour les op\'erateurs 
pseudodiff\'erentiels d'ordre entier.  Sa d\'efinition utilise la trace r\'esiduelle introduite
par Guillemin \cite{Guillemin} et Wodzicki \cite{Wodzicky}.
\begin{definition}
Soit $A\in \Psi^{m}(X;E)$ un op\'erateur pseudodiff\'erentiel inversible d'ordre $m\in \bbZ$, alors
son \textbf{d\'eterminant r\'esiduel} est donn\'e par
\[
   \detr(A):= \exp \left( \Trr(\log A) \right)
\]
pourvu que le logarithme $\log A$ de A soit d\'efini, o\`u $\Trr$ est la trace r\'esiduelle de Guillemin et Wodzicki.   
\label{dr.1}\end{definition}
Derri\`ere cette d\'efinition se cachent deux d\'etails analytiques importants.  D'abord, pour 
pouvoir d\'efinir le logarithme de $A$, il faut supposer que $A$ poss\`ede un angle principal $\theta$,
c'est-\`a-dire un angle $\theta$ tel que le symbole principal 
\[
 \sigma_{m}(A)\in \CI(S^{*}X;\hom(E,E))
\]
ne poss\`ede aucune valeur propre contenue dans la coupure spectrale
\[
                R_{\theta}=\{ re^{i\theta} \quad | \quad r\ge 0\}.
\]
Quoique la d\'efinition du logarithme de $A$ d\'epend du choix de l'angle 
principal, il s'av\`ere que le d\'eterminant r\'esiduel quant \`a lui 
ne d\'epend pas de ce choix (voir \cite{Scott}).
On doit aussi invoquer le r\'esultat de Okikiolu \cite{Okikiolu} pour donner un sens \`a la trace r\'esiduelle
de $\log A$, qui est un op\'erateur pseudodiff\'erentiel logarithmique.   

Lorsqu'on se restreint aux op\'erateurs pseudodiff\'erentiels inversibles d'ordre $0$, on peut toutefois
utiliser une version infinit\'esimale de la d\'efinition~\ref{dr.1} qui
contourne ces difficult\'es analytiques.
\begin{definition}[version infinit\'esimale]
Soit $\cU\subset \Go(X)$ un ouvert simplement connexe contenant l'identit\'e, alors pour $A\in \cU$,
le \textbf{d\'eterminant r\'esiduel} est donn\'e par
\[
   \detr(A):= \exp\left( \int_{0}^{1} \Trr\left[ \gamma^{-1}(t) \frac{d\gamma}{dt}(t)\right] dt \right)
\]
o\`u $\gamma:[0,1]\to \cU$ est une application diff\'erentiable telle que $\gamma(0)=\Id$ et 
$\gamma(1)=A$.
\label{dr.2}\end{definition}
\begin{lemme}
Le d\'eterminant r\'esiduel ne d\'epend pas du choix de l'application diff\'erentiable $\gamma$.  De plus, pour $A$ et $B$ dans $\cU$ suffisamment pr\`es de l'identit\'e, on a que
\[
      \detr(A\circ B) = \detr(A)\detr(B).
\]
\label{dr.3}\end{lemme}
\begin{proof}
Comme on suppose que $\cU$ est simplement connexe, par le th\'eor\`eme de Stokes, il suffit de
v\'erifier que la 1-forme $\Trr[A^{-1}dA]$ d\'efinie sur $\Go(X)$ est ferm\'ee pour voir que la d\'efinition ne d\'epend pas du choix de l'application
$\gamma$.  Or,
\begin{equation*}
\begin{aligned}
d \Trr(A^{-1}dA) &= \Trr( d(A^{-1}dA)) \\
                 &= \Trr(-A^{-1}dA A^{-1}dA) \\
                 &= -\frac{1}{2}\Trr( [A^{-1}dA, A^{-1}dA]) \\
                 &=0, 
\end{aligned}
\end{equation*} 
la derni\`ere \'egalit\'e d\'ecoulant du fait que la trace r\'esiduelle est vraiment une
trace, \`a savoir qu'elle donne z\'ero lorsqu'\'evalu\'ee sur un commutateur.  Cette propri\'et\'e
permet aussi de montrer que le d\'eterminant r\'esiduel est multiplicatif.  Si $A,B\in \cU$ sont 
suffisamment pr\`es de l'identit\'e, alors il existe des applications diff\'erentiables 
$\gamma,\beta:[0,1]\to \cU$ avec $\gamma(0)=\beta(0)=\Id$, $\gamma(1)=A$ et $\beta(1)=B$ telles que $\gamma\beta$ prenne aussi valeur dans $\cU$.  On a alors que
\begin{equation*}
\begin{aligned}
\log \detr(AB)&= \int_{0}^{1} \Trr\left[ (\gamma\beta)^{-1} \frac{d (\gamma\beta)}{dt}\right] dt \\
              &= \int_{0}^{1} \Trr\left[ \beta^{-1}\gamma^{-1} \left(\frac{d\gamma}{dt}\beta +\gamma                       \frac{d\beta}{dt}\right) \right]dt \\
              &= \int_{0}^{1} \Trr\left[ \gamma^{-1}\frac{d\gamma}{dt}\right] dt + 
                   \int_{0}^{1} \Trr\left[ \beta^{-1}\frac{d\beta}{dt}\right] dt \\
              &= \Trr(A)+ \Trr(B),     
\end{aligned}
\end{equation*}
d'o\`u l'on d\'eduit que $\detr(AB)=\detr(A)\detr(B)$.
\end{proof}

Pour $\theta\in (0,2\pi)$, consid\'erons l'ouvert
\[
     \cU_{\theta}:=\{ A\in \Go(X;E)\quad | \quad \theta \;
\mbox{est un angle principal pour} \, A\}\subset
              \Go(X;E).  
\]  
\begin{lemme}
Les d\'efinitions \ref{dr.1} et \ref{dr.2} 
sont \'equivalentes sur un voisinage de l'identit\'e dans $\cU_{\theta}$.  
\label{dede.1}\end{lemme}
\begin{proof}
Par le lemme~\ref{dr.3}, en choisissant notre voisinage suffisamment petit, on a un determinant multiplicatif dans les deux cas. Il suffit alors de 
v\'erifier que la diff\'erentielle de leur logarithme sur le plan tangent
\`a l'identit\'e est la m\^eme.  Pour la d\'efinition~\ref{dr.2}, on voit
directement que 
\[
               \left. d\log \detr\right|_{\Id} = \Trr .
\] 
Pour calculer la diff\'erentielle du logarithme du d\'eterminant dans le 
cas de la d\'efinition \ref{dr.1}, choisissons le voisinage $\mathcal{U}
\subset\cU_{\theta}$ de l'identit\'e suffisamment petit de sorte qu'on ait
pour $A\in \cU$
\[
        \log_{\theta} A=   \sum_{n=1}^{\infty} (-1)^{n-1} \frac{(\gamma-1)^{n}}{n}.
\] 
Si $\gamma: [0,1]\to \cU$ est une application de classe $\CI$
telle que $\gamma(0)=\Id$ et $\gamma(1)=A$, on aura alors que
\begin{equation*}
\begin{aligned}
\Trr \left( \frac{d}{dt} \log_{\theta} \gamma \right) &=
    \Trr\left(  \sum_{n=1}^{\infty} (-1)^{n-1} (\gamma-1)^{n-1} \frac{d\gamma}{dt} 
    \right)  \\
    &= \Trr \left( \gamma^{-1} \frac{d\gamma}{dt} \right),
\end{aligned}
\end{equation*}
d'o\`u l'on d\'eduit que 
\[
     \Trr(\log_{\theta}A)= \Trr\left( \int_{0}^{1} \gamma^{-1}\frac{d\gamma}{dt} dt 
     \right)= \int_{0}^{1} \Trr\left( \gamma^{-1}\frac{d\gamma}{dt} 
     \right) dt,
\]
ce qui montre que la diff\'erentielle du logarithme du d\'eterminant
est aussi donn\'ee par $\Trr$ lorsqu'on utilise la d\'efinition~\ref{dr.1}. 
\end{proof}

En quelque sorte, la version
infinit\'esimale de la d\'efinition du d\'eterminant r\'esiduel remplace la condition de l'existence
d'un angle principal par une condition topologique sur le domaine de 
d\'efinition.  De ce point de vue, on est amen\'e 
\`a se poser la question suivante.
\begin{question}
Est-il possible d'\'etendre la d\'efinition du d\'eterminant r\'esiduel \`a toute la composante
connexe $\Goi(X)$ de l'identit\'e dans $\Go(X)$ par
\begin{equation}
     \detr(A)= \exp \left[ \int_{0}^{1} \Trr\left( \gamma^{-1}\frac{d\gamma}{dt}\right) dt \right]
\label{dr.7}\end{equation}
o\`u $\gamma$ est une application diff\'erentiable telle que $\gamma(0)=\Id$, $\gamma(1)=A$?
\label{dr.4}\end{question}
Comme le lecteur l'aura devin\'e, cette question est purement topologique.  Il suffit de v\'erifier
que cette d\'efinition du d\'eterminant r\'esiduel ne d\'epend pas du choix de l'application 
diff\'erentiable $\gamma$.  On aura une telle ind\'ependance de choix si et seulement si
l'homomorphisme de groupe
\begin{equation}
   \begin{array}{lccl}\Ar: & \pi_{1}(\Go(X)) &\to & \bbC \\
                           & [\gamma] & \mapsto &  \int_{\bbS^{1}} 
                           \Trr\left[\gamma^{-1}\frac{d\gamma}{dt}\right] dt 
   \end{array}                        
\label{dr.5}\end{equation}
prend seulement valeur dans $2\pi i \bbZ$.  En fait, on va montrer
que l'homomorphisme de groupe \eqref{dr.5} est toujours trivial.  L'id\'ee centrale de l'argument que
l'on va pr\'esenter a \'et\'e sugg\'er\'ee  \`a l'auteur par Sergiu Moroianu (voir aussi le paragraphe 8 de \cite{Moroianu} pour une situation similaire) .  Via l'identification $\nu:\pi_{1}(\Go(X))\to K^{0}(S^{*}X)$, on peut voir l'homomorphisme
de groupe $\Ar$ comme une application
\[
                  \Ar: K^{0}(S^{*}X)\to \bbC.
\]
Or, par le biais de la suite exacte \`a six termes 
\begin{equation}
\xymatrix{K_{c}(T^{*}X)\ar[r]&
 K^{0}(X)\ar[r]^{\pi^{*}}&
 K^{0}(S^{*}X)\ar[d]^{\delta}\\
 K^{1}(S^{*}X)\ar[u]&
 K^1(X)\ar[l]&
 K_{c}^{1}(T^{*}X)\ar[l]}
\label{dr.6}\end{equation}
associ\'ee \`a la paire d'espaces $(\overline{T^{*}X},S^{*}X)$, on a une 
inclusion
\begin{equation}
      \pi^{*}(K^{0}(X))\subset  K^{0}(S^{*}X).
\label{dr.8}\end{equation}
\begin{lemme}
Pour tout $[\gamma]\in \pi_{1}(\Go(X))$ tel que $\nu([\gamma])\in 
\pi^{*}(K^{0}(X))$, on a $\Ar([\gamma])=0$. 
\label{dr.13}\end{lemme}
\begin{proof}
Comme $G^{-\infty}(\bbS^{1})$ est un espace classifiant pour la $K$-th\'eorie
impaire, on a l'identification
\[
    K^{0}(X)\cong K^{-2}(X)\cong [(S^{1}(X^{+}),\pt); (G^{-\infty}(\bbS^{1},\Id)]
\]
et l'application $\pi^{*}:K^{-2}(X)\to K^{-2}(S^{*}X)$ est alors induite par 
l'image r\'eciproque associ\'ee \`a la projection
\[
        \pi: S^{1}(S^{*}X^{+})\to S^{1}(X^{+})
\]
o\`u $X^{+}=X\cup \pt$ est l'union disjointe de $X$ avec un point.  En faisant appel au 
th\'eor\`eme~\ref{gh.3} cela 
montre que lorsque $\nu([\gamma])\in \pi^{*}(K^{0}(X))$, on peut choisir
$\gamma\in \CI([0,1]; \Go(X))$ repr\'esentant $[\gamma]\in \pi_{1}(\Go(X))$
de sorte que pour tout $t\in [0,1]$, $\gamma(t)\in \Go(X)$ soit simplement
donn\'e par un isomorphisme de fibr\'es vectoriels.  En ce cas, le terme 
d'ordre $-\dim X$ du symbole total (full symbol en anglais) de $\gamma^{-1}
\frac{d\gamma}{dt}$ est nul, ce qui signifie que 
\[
           \Trr \left( \gamma^{-1}\frac{d\gamma}{dt}\right)=0
\]
pour tout $t\in [0,1]$ \'etant donn\'ee la formule bien connue
\begin{equation}
       \Trr(A)= \frac{1}{(2\pi)^{n}} \int_{X}\int_{|\xi|=1} \operatorname{tr}
\sigma(A)_{-n}(x,\xi) dS(\xi)dx, \quad n=\dim X.
\label{star.1}\end{equation}
exprimant le trace r\'esiduelle en termes de la partie d'ordre $-n$ du symbole
total.  La d\'efinition de ce terme d'ordre $-n$ n'est pas naturelle
et d\'epend d'un choix de quantification, mais il s'av\`ere que la formule \eqref{star.1}
est ind\'ependante d'un tel choix.

On obtient donc
\[
              \Ar([\gamma])= \int_{0}^{1} \Trr\left( \gamma^{-1}
\frac{d\gamma}{dt}\right) dt=0.
\]
\end{proof}
Revenant \`a la suite exacte \`a six termes \eqref{dr.6}, rappelons
que l'homomorphisme de bord $\delta: K^{0}(S^{*}X)\to K^{1}_{c}(T^{*}X)$ est
toujours surjectif (voir par exemple p.81 dans \cite{APS3}).  Toujours selon \cite{APS3} (p.81), on 
peut repr\'esenter un \'el\'ement de $K^{-1}_{c}(T^{*}X)$ par un lacet
\begin{equation}
\begin{aligned}
   \sigma_{t} &= \Id \cos t + i\sigma \sin t , \quad 0\le t\le \pi, \\
              &= \Id( \cos t + i\sin t), \quad \pi \le t \le 2\pi,
\end{aligned}
\label{dr.9}\end{equation}
o\`u $\sigma\in \CI(S^{*}X; G^{-\infty}(\bbS^{1}))$ est un symbole auto-adjoint.
On peut aussi supposer via l'identification \eqref{stab.4} que le symbole $\sigma$
est en fait un \'el\'ement de 
\[
\CI(S^{*}X;\Hom(\pi^{*}E,\pi^{*}E)) 
\]
pour un certain
fibr\'e vectoriel complexe $E\to X$ o\`u $\pi: S^{*}X\to X$ est la projection de fibr\'e.
Soit $A\in \Psi^{1}(X;E)$ un op\'erateur auto-adjoint inversible d'ordre $1$ dont le symbole est donn\'e par $\sigma$.  Le spectre de $A$ est donc discret et contenu
dans $\bbR$.  Comme l'op\'erateur $A$ est inversible, les projections $P_{+}$ et
$P_{-}$ sur les espaces propres positifs et negatifs sont des op\'erateurs 
pseudodiff\'erentiels d'ordre z\'ero, comme on peut le voir en les d\'efinissant 
\`a la mani\`ere de Seeley \cite{Seeley} par des int\'egrales de contour dans
le plan complexe utilisant le r\'esolvant de $A$.  Ainsi,
\[
       |A|:= P_{+}A - P_{-}A\in \Psi^{1}(X;E)
\]
est un op\'erateur pseudodiff\'erentiel inversible d'ordre $1$.  Son symbole est donn\'e
par $|\sigma|= \sqrt{\sigma^{2}}$, un symbole auto-adjoint dont les valeurs propres sont toutes
positives pour tout $x\in S^{*}X$.  Le symbole $|\sigma|$ est donc homotope au symbole 
identit\'e dans l'espace des symboles auto-adjoints inversibles.  Consid\'erons donc
l'op\'erateur
\[
     A_{0}:= |A|^{-1}A\in \Psi^{0}(X;E).
\]   
Clairement, on a que $A_{0}^{2}= \Id$ et le spectre de $A_{0}$ est contenu dans 
$\{-1,1\}$.  En fait, si $f\in \CI(X;E)$ est une section propre de $A$,
\[
                 Af= \lambda f, \quad \lambda \in \bbR\setminus\{0\},
\]
alors $f$ est aussi une section propre de $A_{0}$ avec 
\[
        A_{0}f= \frac{\lambda}{|\lambda|} f.
\]
Quant \`a lui, le symbole $\sigma(A_{0})$ de $A_{0}$ est homotope au symbole $\sigma$ dans
l'espace des symboles auto-adjoints inversibles.  L'\'el\'ement de $K^{-1}_{c}(T^{*}X)$ 
repr\'esent\'e par le lacet \eqref{dr.9} peut donc aussi 
\^etre repr\'esent\'e par $\sigma(\gamma_{t})$ o\`u
\begin{equation}
 \begin{aligned}
   \gamma_{t} &= \Id \cos t + iA_{0} \sin t , \quad 0\le t\le \pi, \\
              &= \Id( \cos t + i\sin t), \quad \pi \le t \le 2\pi,
\end{aligned}
\label{dr.26}\end{equation}

\begin{theorem}
Pour toute vari\'et\'e compacte sans bord $X$ de classe $\CI$, l'homomorphisme de bord
\[
      \Ar: \pi_{1}(\Go(X))\to \bbC
\]
est trivial.  La d\'efinition~\eqref{dr.2} du d\'eterminant r\'esiduel peut donc toujours \^etre 
\'etendue globalement \`a toute la composante connexe $\Goi(X)$ de l'identit\'e dans $\Go(X)$ pour 
donner lieu \`a un d\'eterminant multiplicatif.
\label{dr.11}\end{theorem}
\begin{proof}
Soit $[\gamma]\in \pi_{1}(\Go(X))$.  On veut montrer que $\Ar([\gamma])=0$.  Par le lemme~\ref{dr.13}, 
le th\'eor\`eme~\ref{gh.3} et la discussion qui pr\'ec\`ede, on peut supposer que 
$[\gamma]$ est repr\'esent\'e par un lacet de la forme \eqref{dr.26}.  Pour $t\in [\pi,2\pi]$,
\begin{equation*}
\begin{aligned}
   \gamma^{-1}(t) \frac{d\gamma}{dt} &= \Id (\cos t -i \sin t)( -\sin t+ i\cos t) \\
                               &= i(\Id)
\end{aligned}
\end{equation*}
est un multiple de l'identit\'e, donc 
\[
         \Trr\left( \gamma^{-1}(t) \frac{d\gamma}{dt} \right)=0
\]
dans ce cas, puisque comme dans la preuve du lemme~\ref{dr.13}, la partie d'ordre $-n$ du symbole 
total est nulle.  Ainsi, on a que 
\[
     \Ar([\gamma])= \int_{0}^{\pi} \Trr\left( \gamma^{-1} \frac{d\gamma}{dt}\right) dt.
\]
Par d\'efinition de la trace r\'esiduelle, cette quantit\'e est donn\'ee par le r\'esidu
en $s=0$ de l'extension m\'eromorphe \`a tout le plan complexe de la fonction holomorphe
\[
   s\mapsto \xi(s):= \int_{0}^{\pi} \Tr\left( |A|^{-s} \gamma^{-1} \frac{d\gamma}{dt} \right) dt,
   \quad s\in \bbC, \quad \Re s>>0.
\]
Comme $A_{0}^{-1}= A_{0}$, on a pour $t\in [0,\pi]$ que
\begin{equation*}
\begin{aligned}
\gamma^{-1} \frac{d\gamma}{dt} &= (\Id \cos t + iA_{0}\sin t)^{-1}( -\Id\sin t+ iA_{0} \cos t) \\
                               &= (\Id\cos t -i A_{0}\sin t)( -\Id \sin t + iA_{0} \cos t) \\
                               &= \Id(-\cos t \sin t + \cos t \sin t) + iA_{0}(\cos^{2} t + \sin^{2} t)\\
                               &= i A_{0}.
\end{aligned}
\end{equation*}
Ainsi, pour $\Re s>>0$, on a 
\begin{equation*}
\begin{aligned}
\xi(s) &= \int_{0}^{\pi} \Tr( |A|^{-s}i A_{0}) dt \\
                               &= i\pi \Tr(|A|^{-s}A_{0})= i\pi \Tr(|A|^{-s-1}A) \\
                               &= i\pi \sum_{\lambda\in \spec(A)} |\lambda|^{-s-1}\lambda\\
                               &= i \pi \eta(A,s)
\end{aligned}
\end{equation*}
o\`u $\eta(A,s)$ est  la fonctionnelle de Atiyah Patodi et Singer \cite{APS1} associ\'ee \`a l'op\'erateur auto-adjoint
inversible $A\in \Psi^{1}(X;E)$.  D'apr\`es les r\'esultats de Gilkey \cite{Gilkey} et Wodzicki \cite{Wodzicky_eta},
\cite{Wodzicki_eta2}, la fonctionnelle
$\eta(A,s)$ n'a pas de pole \`a $s=0$, d'o\`u l'on conclut que $\Ar([\gamma])=0$.  Le d\'eterminant
r\'esiduel est donc d\'efini globalement sur $\Goi(X)$ et la seconde partie de la d\'emonstration du lemme~\ref{dr.3} montre alors que c'est un d\'eterminant multiplicatif.
\end{proof}
Plus g\'en\'eralement, pour $n\in \bbZ$ et $E\to X$ un fibr\'e vectoriel complexe, on peut consid\'erer
l'espace 
\[
    G^{n}(X;E):= \{ A\in \Psi^{n}(X;E)\quad | \quad A \; \mbox{est inversible} \}
\]
des op\'erateurs inversibles d'ordre $n$.  La composition d'op\'erateurs induit une structure de groupe
sur
\[
   G^{*}(X;E):= \bigcup_{n\in \bbZ} G^{n}((X;E).
\]
\begin{corollaire}
Soit $P\in \Psi^{1}(X;E)$ un choix d'op\'erateur auto-adjoint inversible dont les valeurs propres 
sont toutes positives.  Alors pour $n\in\bbZ$, le d\'eterminant r\'esiduel admet une d\'efinition
globale sur la composante connexe $G^{n}_{P}(X;E)$ de $P^{n}$ dans $G^{n}(X;E)$.
\label{dr.12}\end{corollaire}
\begin{proof}
Clairement, $P^{n}$ poss\`ede un angle principal.  On peut donc d\'efinir $\detr(P^{n})$ en utilisant
la d\'efinition~\ref{dr.1}.  Pour $A\in G^{n}_{P}(X;E)$ quelconque, on a alors que 
$AP^{-n}\in \Goi(X;E)$ et on pose donc
\[
    \detr(A)= \detr(AP^{-n})\detr(P^{n}).
\]
\end{proof}

\bibliographystyle{amsplain}

\providecommand{\bysame}{\leavevmode\hbox to3em{\hrulefill}\thinspace}
\providecommand{\MR}{\relax\ifhmode\unskip\space\fi MR }
\providecommand{\MRhref}[2]{%
  \href{http://www.ams.org/mathscinet-getitem?mr=#1}{#2}
}
\providecommand{\href}[2]{#2}

\end{document}